\documentclass[compress,final,3p,times,18pt]{elsarticle}
\usepackage{mathrsfs}
\usepackage{amsfonts}
\usepackage{graphics}
\usepackage{amssymb}
\usepackage{amsmath}
\usepackage{subfig}
\usepackage{float}
\usepackage{caption}
\usepackage{dsfont}

\allowdisplaybreaks

\newtheorem{Theo}{Theorem}
\newtheorem{Assu}{Assumption}
\newtheorem{Rem}{Remark}
\newtheorem{Lem}{Lemma}

\newtheorem{Exa}{Example}
\numberwithin{equation}{section}
\numberwithin{Assu}{section}
 \numberwithin{Lem}{section}
 \numberwithin{Defi}{section}
 \numberwithin{Theo}{section}
 \numberwithin{Rem}{section}
  \numberwithin{Coro}{section}
  \numberwithin{Fig}{section}

\journal{}

\begin{document}

\begin{frontmatter}



\title{Convergence and stability of the semi-tamed Milstein method for commutative stochastic differential equations with non-globally Lipschitz continuous coefficients\tnoteref{label1}} \tnotetext[label1]{This work was supported by Natural Science Foundation of Hunan Province (Grant No. 2018JJ3628), National Natural Science Foundation of China (Grant Nos. 12071488, 11971488 and 11901527) and the China Postdoctoral Science Foundation (Grant No. 2020M671087).}
\author{Yulong Liu\fnref{addr1}}
\ead{liuyulong@csu.edu.cn}
\author{Yuanling Niu\corref{cor1}\fnref{addr1}}
\ead{yuanlingniu@csu.edu.cn}
\author{Xiujun Cheng\fnref{addr2,addr3}}
\ead{xiujuncheng@zstu.edu.cn} \cortext[cor1]{Corresponding author}

\address[addr1]{School of Mathematics and Statistics, Central South University, Changsha, 410083, China}
\address[addr2]{College of Science, Zhejiang Sci-Tech University, Hangzhou, 310018,  China}
\address[addr3]{Institute of Natural Sciences, Shanghai Jiao Tong University, Shanghai, 200240, China}

\begin{abstract}
{A new explicit stochastic scheme of order 1 is proposed for solving commutative stochastic differential equations (SDEs) with non-globally Lipschitz continuous coefficients.
The proposed method is a semi-tamed version of Milstein scheme to solve SDEs with the drift coefficient consisting of non-Lipschitz continuous term and globally Lipschitz continuous term.
It is easily implementable and achieves higher strong convergence order. A stability criterion for this method is derived, which shows that the stability condition of the numerical methods and that of the solved equations keep uniform. Compared with some widely used numerical schemes, the proposed method has better performance in inheriting the mean square stability of the exact solution of SDEs.  Numerical experiments are given to illustrate the obtained convergence and stability properties.}
\end{abstract}


\end{frontmatter}

\section{Introduction}
\label{sec:into}
In many scientific fields, such as economics, finance and complex networked systems,
stochastic differential equations (SDEs) are often used to model complex
dynamics \cite{Mao1997,zong2020}. It has become an important issue to develop numerical methods for SDEs since exact solutions of SDEs can rarely be obtained.
In this work, we consider the numerical solution for the SDEs in the It\^{o} sense
\begin{equation}
	\label{eq_1_1}
	dX_t = f(X_t)dt + g(X_t)dW_t, ~~X_0 = \xi,~~ t\in [0,T].
\end{equation}
Here
$f : \mathbb{R}^d \rightarrow \mathbb{R}^d,g = (g_1, g_2,\cdots, g_m): \mathbb{R}^d \rightarrow \mathbb{R}^{d\times m}.$
Assume that $g$ is global Lipschitz continuous and $f$ is globally one-sided Lipschitz continuous. Moreover, the drift coefficient $f$ has the form: $ f(x) = \phi(x)+ \varphi(x)$, where $\phi: \mathbb{R}^d \rightarrow \mathbb{R}^d$ is Lipschitz continuous, and $\varphi: \mathbb{R}^d \rightarrow \mathbb{R}^d$ is non-Lipschitz continuous.
As usual, we assume that $W_t$ is an m-dimensional Brownian motion defined on the complete probability space $(\Omega, \mathcal{F}, \mathbb{P})$ equipped with an increasing filtration $\{\mathcal{F}_t\}_{t\geqq 0}$ satisfying the usual conditions. The initial data $\xi$ is independent of the Brownian motion. Eq. \eqref{eq_1_1} can be interpreted mathematically as the following integral equation
\begin{equation}
	\label{eq_1_2}
	X_t = X_0 + \int_{0}^{t} f(X_s)ds + \sum_{i=1}^{m}\int_{0}^{t} g_i(X_s)dW_s, t\in [0,T],
\end{equation}
where $g_i(x) = (g_{1,i}, g_{2,i},\cdots, g_{d,i})^T$ for $x \in \mathbb{R}^d, i \in \{1,2, \cdots, m\}$
and the second integral is the It\^{o}'s integral.

The work is concerned with the strong approximation of SDEs \eqref{eq_1_1} and its convergence as well as the mean square stability. We establish a numerical approximation called semi-tamed Milstein method $Y^N_n: \Omega \rightarrow \mathbb{R}^d, n \in \{0,1,\cdots,N\}$ with $Y^N_0 = \xi $ on a uniform mesh with stepsize $h = T/N$ defined by $\mathcal{T}^N: 0 = t_0 < t_1 < \cdots <t_N = T, N \in \mathbb{N}$. Its strong convergence order given by
\begin{equation}
	\label{eq_1_3}
	\left(\mathbb{E}\bigg[\sup_{t\in[0,T]}\|X_t - \bar{Y}^N_t\|^p \bigg] \right)^{1/p} \leq C_{p,T}\cdot h,\qquad h\in(0,1], \qquad C_{p,T} \geq 1, \qquad p \in[1,\infty)
\end{equation}
is obtained.
We also give the mean square stability of the numerical approximation, which shows better stability preservation of the underlying method under a stepsize restriction.

The strong approximation of SDEs \eqref{eq_1_1} has been studied extensively. The classic explict Euler-Maruyama (EM) method, proposed by Maruyama \cite{maruyama1955}, investigated in \cite{kloeden1992} and \cite{milstein1995}, is strongly convergent with order one-half if the coefficients $f, g$ satisfy the global Lipschitz continuous condition. Yuan and Mao \cite{Yuan-Mao2008} obtained the convergence rate of explicit EM method for locally Lipschitz continuous coefficients. Unfortunately, Hutzenthaler et al. \cite{hutzenthaler2009, hutzenthaler2011} showed that the explicit EM approximation for SDEs with superlinearly growing drift fails to converge strongly to the exact solution. Although the backward Euler scheme can hold the convergence in such a situation, it needs extra computational effort since it is implicit. To solve this problem, the ``tamed" methods such as tamed Euler scheme \cite{hutzenthaler2012} and tamed Milstein scheme \cite{wang2013} were proposed, which are explicit schemes and strongly convergent to exact solution of SDEs with non-Lipschitz continuous drift coefficients. Sabanis \cite{sabanis2016} proposed a ``tamed" version explicit method for SDEs with superlinearly growing drift and diffusion coefficients and studied its convergence in probability and in $L^p$ sense. Moreover, the ``truncated" method such as the truncated Euler-Maruyama method in \cite{mao2015} and the truncated Milstein method in \cite{guo2018, zhang2019} can also well resolve the super-linear growth of the coefficients.

The stability analysis of the numerical methods for SDEs plays an important role in computation and simulation. For exact solutions of SDEs, Mao \cite{Mao1997} gave the concepts of T-stability, $p$-th moments stability et al. On the other hand, for the numerical solutions of nonlinear SDEs, the stability, especially the mean square stability has been widely studied for different numerical methods, e.g. see \cite{zong2014,Gan2014,Yao2018}. But there is few stability analysis of tamed version schemes.

Inspired by literature \cite{wang2013} and \cite{zong2014}, we derive a semi-tamed Milstein method which is a non-Lipschitz term tamed approximation for SDEs \eqref {eq_1_1}
\begin{align}
	\label{eq_1_4}
	Y^N_{n+1} = Y^N_n& + \phi(Y^N_n)h + \frac{\varphi(Y^N_n)h}{1 + \|\varphi(Y^N_n)\|h} + g(Y^N_n)\Delta W_n  \notag\\
	 & + \frac{1}{2}\displaystyle\sum_{j_1,j_2=1}^{m}L^{j_1}g_{j_2}(Y^N_n)(\Delta W^{j_1}_n\Delta W^{j_2}_n - \delta_{j_1,j_2}h),
\end{align}
where $f(x) = \phi(x) + \varphi(x)$ and $L^{j_1} = \displaystyle\sum_{k=1}^{d}g_{k,j_1}\frac{\partial}{\partial x^k}$.  $\Delta W_n = W_{(n+1)h} - W_{nh}$ is the m-dimensional Brownian increment, $\Delta W^{j_1}_n, \Delta W^{j_2}_n$ is the mutual independent Brownian increment for $j_1,j_2 \in \{1,2,\cdots,m\}$.
 Here, to propose a more
computationally efficient Milstein-type method in this work, we only consider SDEs with diffusion matrix $g$ fulfills the so-called commutativity condition
\begin{equation}
	\label{eq_1_5}
	L^{j_1}g_{k,j_2} = L^{j_2}g_{k,j_1}.
\end{equation}
The underlying explicit method \eqref{eq_1_4} is shown to have higher strong convergence order $1$ and better mean square stability compared with some widely used numerical schemes.

The paper is organized as follows. Uniform boundness of $p$th moments is presented in
Section 2. Strong convergence order of the semi-tamed Milstein method is shown in Section 3. Exponential mean square stability of the method is given in Section 4. Numerical
results are reported confirming convergence properties and comparing stability properties in Section 5.

\section{Uniform boundedness of $p$th moments}
\label{sec:uniform boundedbess}
We use the following notations throughout this paper. Let $N \in \mathbb{N}$ be the step number of the uniform mesh defined in the previous section. $T \in (0,\infty)$ is a fixed real number and $h = T/N$ is the stepsize. Moreover, we define $\left\|x\right\| = (|x_1|^2 + \cdots + |x_k|^2)^{1/2}$, $\left<x,y\right> = x_1y_1 + \cdots + x_ky_k$ for all $x = (x_1, x_2, \cdots, x_k), y = (y_1, \cdots, y_k) \in \mathbb{R}^k, k \in \mathbb{N}$, and $\left\|A\right\| = \sup_{x\in\mathbb{R}^l,||x||\leq 1}\left\|Ax\right\|$ for all $A \in \mathbb{R}^{k\times l}, k,l \in \mathbb{N}$.
\begin{Assu}
	\label{assumption2_1}
Let $\varphi(x)$ and $g_i(x), i \in \{1,2,\cdots,m\}$ be continuously differentiable functions and there exist positive constants $K \geq 1$ and $c \geq 1$, such that for all $ x, y \in \mathbb{R}^d$, the following conditions hold:
\begin{align}
\label{eq_2_1}\left<x-y,f(x)-f(y)\right> \leq K\left\|x-y\right\|^2,\qquad\qquad\qquad\qquad   \\
\label{eq_2_2}\left\|g(x)-g(y)\right\| \leq K\left\|x-y\right\|,\qquad\qquad\qquad\qquad \\
\label{eq_2_3}\left\|\phi(x)-\phi(y)\right\| \leq K\left\|x-y\right\|,\qquad\qquad\qquad\qquad           \\
\label{eq_2_4}\left\|\varphi'(x)\right\| \leq K\left\|x\right\|^c,\qquad\qquad\qquad\qquad\qquad\qquad   \\
\label{eq_2_5}\left\|L^{j_1}g_{j_2}(x)-L^{j_1}g_{j_2}(y)\right\| \leq K\left\|x-y\right\|, j_1, j_2 \in \{1,2,\cdots,m\}.
\end{align}
\end{Assu}

\begin{Rem}
Inequality \eqref{eq_2_1} in Assumption \ref{assumption2_1}, named as one-sided Lipschitz condition, has been widely used, for instance, in \cite{Burrage1979, gyongy1998}. Moreover, the drift coefficient of many stochastic models, such as stochastic Ginzburg-Landau equation \cite{kloeden1992, hutzenthaler2011b}, stochastic Lorenz equation \cite{hutzenthaler2011b, schmallfus1997} and volatility process \cite{kloeden1992, hutzenthaler2011b} has the form: $f(x) = \phi(x) + \varphi(x)$, where the part $\phi$ satisfies condition \eqref{eq_2_3} and the part $\varphi$ satisfies condition \eqref{eq_2_4} respectively.
\end{Rem}

We introduce the dominating stochastic process $D_{n}$ and appropriate sub-events $\Omega_{n}$ to prove the uniform boundedness of pth moments of our numerical solution:
\begin{align}
    \label{eq_2_7}	D_n &:= (\lambda + \|\xi\|)\exp(\lambda + \displaystyle\sup_{0\leq u\leq n} \displaystyle\sum_{k=u}^{n-1}[\lambda\|\Delta W_k\|^2 + \alpha_k]),
    \\
    \label{eq_2_6}\Omega_n &:= \left\{\omega \in \Omega| \displaystyle\sup_{0\leq k\leq n-1} D_k(\omega)\leq N^{1/2c}, \displaystyle\sup_{0\leq k\leq n-1}\|\Delta W_k\| \leq 1 \right\},
\end{align}
where
\begin{align}
	\lambda = \bigg( &1 + 4TK + 2T(\|\phi(0)\| + \|\varphi(0)\|) + 2K + 2\|g(0)\| \notag\\
                     &  + m^2(T+1)\left(K + \displaystyle\max_{1\leq j_1,j_2\leq m}\|L^{j_1}g_{j_2}(0)\|\right) \bigg)^2,
\end{align}
and
\begin{align}
	\label{eq_2_8}
	\alpha_n &= \mathds{1}_{\{\|Y^N_n\|\geq 1\}} \left\langle\frac{Y^N_n}{\|Y^N_n\|}, \frac{g(Y^N_n)}{\|Y^N_n\|}\Delta W_n\right\rangle\notag  \\
			 &+ \mathds{1}_{\{\|Y^N_n\|\geq 1\}} \left\langle\frac{Y^N_n}{\|Y^N_n\|}, \displaystyle\sum_{j_1,j_2=1}^{m} \frac{L^{j_1}g_{j_2}(Y^N_n)}{2\|Y^N_n\|}\left(\Delta W_{n}^{j_1}\Delta W_{n}^{j_2} - \delta_{j_1,j_2}h \right)\right\rangle.
\end{align}
For convenience, we also make the following notation
\begin{equation}
	\tilde{\varphi}(x) = \frac{\varphi(x)}{1 + h\|\varphi(x)\|}.
\end{equation}

\begin{Lem}\label{lemma2_1}
	 $Y^N_n, D_n$ and $\Omega_n$ are defined in \eqref{eq_1_4}, \eqref{eq_2_7} and \eqref{eq_2_6}. Then we have
	\begin{equation}
		\label{eq_2_9}
		\mathds{1}_{\Omega_n}\|Y^N_n\| \leq D_n,\qquad\text{for all } n \in \{0,1,\cdots, N\}.
	\end{equation}
\end{Lem}

\begin{proof}
Firstly, by the definition of $\Omega_{n}$, we can obtain that $\|\Delta W_n\| \leq 1$ on $\Omega_{n+1}$ for all $0\leq n\leq N-1, N \in \mathbb{N}.$ Then the Lipschitz continuity condition on $g, \phi$ and $L^{j_1}g_{j_2}, j_1,j_2 \in \{1,2,\cdots,m\},$ and the polynomial growth bound condition on $\varphi'$ imply that
\begin{align}
	\label{eq_2_10}
\|Y^N_{n+1}\| &\leq \|Y^N_n\| + h\|\phi(Y^N_n)\| + h\|\tilde{\varphi}(Y^N_n)\| + \|g(Y^N_n)\|\|\Delta W_n\|  \notag\\
			&\qquad + \frac{1}{2}\displaystyle\sum_{j_1,j_2=1}^{m}\|L^{j_1}g_{j_2}(Y^N_n)\||\Delta W^{j_1}_n\Delta W^{j_2}_n - \delta_{j_1,j_2}h|  \notag\\
			&\leq 1 + h\|\phi(Y^N_n) - \phi(0)\| + h\|\phi(0)\| + h\|\varphi(Y^N_n) - \varphi(0)\| \notag\\
			&\qquad + h\|\varphi(0)\| + \|g(0)\|\|\Delta W_n\| + \|g(Y^N_n) - g(0)\|\|\Delta W_n\| \notag\\
			&\qquad + \frac{1}{2}\displaystyle\sum_{j_1,j_2=1}^{m}\left(\|L^{j_1}g_{j_2}(Y^N_n) - L^{j_1}g_{j_2}(0)\| + \|L^{j_1}g_{j_2}(0)\| \right) \times |\Delta W^{j_1}_n\Delta W^{j_2}_n - \delta_{j_1,j_2}h| \notag \\	
			&\leq 1 + hK\|Y^N_n\| +h\|\phi(0)\| + hK\|Y^N_n\|^{c+1} + h\|\varphi(0)\| + K\|Y^N_n\|\|\Delta W_n\| + \|g(0)\|\|\Delta W_n\|  \notag\\
			&\qquad + \frac{1}{2}\displaystyle\sum_{j_1,j_2=1}^{m}\left(K\|Y^N_n\| + \|L^{j_1}g_{j_2}(0)\| \right) \times |\Delta W^{j_1}_n\Delta W^{j_2}_n - \delta_{j_1,j_2}h| \notag\\
			&\leq 1 + hK(1 + \|Y^N_n\|^c)\|Y^N_n\| + h(\|\phi(0)\| + \|\varphi(0)\|) + K + \|g(0)\|  \notag\\
			&\qquad + \frac{m}{2}(m + T)\left(K + \displaystyle\max_{1\leq j_1,j_2\leq m}\|L^{j_1}g_{j_2}(0)\|\right)  \notag\\
			&\leq 1 + 2TK + T(\|\phi(0)\| + \|\varphi(0)\|) + K + \|g(0)\|   \notag\\
			&\qquad + \frac{m}{2}(m + T)\left(K + \displaystyle\max_{1\leq j_1,j_2\leq m}\|L^{j_1}g_{j_2}(0)\|\right)  \notag\\
			&\leq \lambda,
\end{align}
on $\Omega_{n+1} \cap \{\omega \in \Omega|\|Y^N_n(\omega)\|\leq 1\}$ for all $0\leq n\leq N-1$.
Moreover, using Cauchy-Schwarz inequality and the inequality $a\cdot b \leq \frac{a^2}{2} + \frac{b^2}{2}$ for all $a,b \in \mathbb{R}$, we can obtain that
\begin{align}
	\label{eq_2_11}
	\|Y^N_{n+1}\|^2 &= \|Y^N_n + \phi(Y^N_n)h + \tilde{\varphi}(Y^N_n)h + g(Y^N_n)\Delta W_n + A_n\|^2 \notag \\
				  &= \|Y^N_n\|^2 + \|\phi(Y^N_n)\|^2h^2 + \|\tilde{\varphi}(Y^N_n)\|^2h^2 + \|g(Y^N_n)\Delta W_n\|^2 + \|A_n\|^2  \notag\\
				  &\quad + 2h\langle Y^N_n, \phi(Y^N_n)\rangle + 2h\langle Y^N_n, \tilde{\varphi}(Y^N_n) \rangle + 2\langle Y^N_n, g(Y^N_n)\Delta W_n \rangle + 2\langle Y^N_n, A_n \rangle   \notag\\
				  &\quad + 2\langle \phi(Y^N_n)h, \tilde{\varphi}(x)h \rangle + 2\langle \phi(Y^N_n)h, g(Y^N_n)\Delta W_n \rangle + 2\langle \phi(Y^N_n)h, A_n \rangle  \notag\\
				  &\quad + 2\langle \tilde{\varphi}(x)h, g(Y^N_n)\Delta W_n \rangle + 2\langle \tilde{\varphi}(x)h, A_n \rangle + 2\langle g(Y^N_n)\Delta W_n, A_n \rangle  \notag\\
				  &\leq \|Y^N_n\|^2 + 4\|\phi(Y^N_n)\|^2h^2 + 4\|\varphi(Y^N_n)\|^2h^2 + 4\|g(Y^N_n)\Delta W_n\|^2 + 4\|A_n\|^2  \notag\\
				  &\quad + 2h\langle Y^N_n, \phi(Y^N_n)\rangle + 2h\langle Y^N_n, \tilde{\varphi}(Y^N_n) \rangle + 2\langle Y^N_n, g(Y^N_n)\Delta W_n \rangle + 2\langle Y^N_n, A_n \rangle,
\end{align}
on $\Omega_{n}$ and $0\leq n\leq N-1.$ Here we denote
\begin{equation}
	\label{eq_2_12}
	A_n = \frac{1}{2}\displaystyle\sum_{j_1,j_2=1}^{m}L^{j_1}g_{j_2}(Y^N_n)(\Delta W^{j_1}_n\Delta W^{j_2}_n - \delta_{j_1,j_2}h).
\end{equation}
The globally Lipschitz condition on $g$ and $L^{j_1}g_{j_2}$ and the one-sided Lipschitz condition on $f$ imply that for $\|x\|\geq 1$
\begin{align}
	\label{eq_2_13}
	\|g(x)\|^2 &\leq (\|g(x) - g(0)\| + \|g(0)\|)^2  \notag\\
	                &\leq (K\|x\| + \|g(0)\|)^2 \notag\\
	                &\leq (K + \|g(0)\|)^2\|x\|^2,
\end{align}
\begin{align}
	\label{eq_2_14}
	\|L^{j_1}g_{j_2}(x)\| &\leq (\|L^{j_1}g_{j_2}(x) - L^{j_1}g_{j_2}(0)\| + \|L^{j_1}g_{j_2}(0)\|)  \notag\\
	 						   &\leq (K\|x\| + \|L^{j_1}g_{j_2}(0)\|)^2  \notag\\
	 						   &\leq (K + \|L^{j_1}g_{j_2}(0)\|)^2\|x\|^2,
\end{align}
\begin{align}
	\label{eq_2_15}
	\langle x, f(x)\rangle &= \langle x, f(x) - f(0)\rangle + \langle x, f(0)\rangle  \notag\\
							 &\leq K\|x\|^2 + \|x\|\cdot\|f(0)\|  \notag\\
							 &\leq (K + \|f(0)\|)\|x\|^2.
\end{align}
Additionally, the globally Lipschitz continuity condition on $\phi$ gives that on $\|x\|\geq 1$
\begin{align}
	\label{eq_2_16}
	\|\phi(x)\|^2 &\leq (\|\phi(x) - \phi(0)\| + \|\phi(0)\|)^2  \notag\\
		       &\leq (K\|x\| + \|\phi(0)\|)^2  \notag\\
		       &\leq (K + \|\phi(0)\|)^2\|x\|^2,
\end{align}
\begin{align}
	\label{eq_2_17}
	|\langle x, \phi(x)\rangle| &\leq \|x\|\|\phi(x) - \phi(0)\| - \|x\|\|\phi(0)\|  \notag\\
							 &\leq K\|x\|^2 + \|\phi(0)\|\|x\|  \notag\\
							 &\leq (K + \|\phi(0)\|)\|x\|^2.
\end{align}
From the polynomial growth bound condition on $\varphi'$, we obtain that on $1 \leq \|x\| \leq N^{1/2c}$
\begin{align}
	\label{eq_2_18}
	\|\varphi(x)\|^2 &\leq (\|\varphi(x) - \varphi(0)\| + \|\varphi(0)\|)^2  \notag\\
			   &\leq (K\|x\|^{c+1} + \|\varphi(0)\|)^2  \notag\\
			   &\leq (K + \|\varphi(0)\|)^2\|x\|^{2(c+1)}  \notag\\
			   &\leq N(K + \|\varphi(0)\|)^2\|x\|^2.
\end{align}
Combining \eqref{eq_2_13}-\eqref{eq_2_18} and the fact that
\begin{equation*}
	 2h\langle Y^N_n, \tilde{\varphi}(Y^N_n) \rangle = \frac{2h\langle Y^N_n, \varphi(Y^N_n)\rangle}{1 + \|\varphi(x)\|h} = \frac{2h\langle Y^N_n, f(Y^N_n)\rangle}{1 + \|\varphi(x)\|h} - \frac{2h\langle Y^N_n, \phi(Y^N_n)\rangle}{1 + \|\varphi(x)\|h},
\end{equation*}
we get from \eqref{eq_2_11} that on $\{\omega \in \Omega|1\leq \|Y^N_n(\omega)\|\leq N^{1/2c}\}$
\begin{align}
	\label{eq_2_19}
	\|Y^N_{n+1}\|^2 &\leq \|Y^N_n\|^2 + 4h^2(K + \|\phi(0)\|)^2\|Y^N_n\|^2 + 4h^2N(K + \|\varphi(0)\|)^2\|Y^N_n\|^2  \notag\\
				&\quad + 4(K + \|g(0)\|)^2\|\Delta W_n\|^2\|Y^N_n\|^2 + 4h(K + \|\phi(0)\|)\|Y^N_n\|^2  \notag\\
				&\quad + 2h(K + \|f(0)\|)\|Y^N_n\|^2 + 2\langle Y^N_n, g(Y^N_n)\Delta W_n \rangle + 2\langle Y^N_n, A_n \rangle + 4\|A_n\|^2.
\end{align}
For the last term $\|A_n\|^2$, using the inequality $\left(\displaystyle\sum_{i=1}^{m}a_i\right)^2 \leq m\left(\displaystyle\sum_{i=1}^{m}a_i^2\right)$, we can obtain the estimate on $\omega \in \Omega_{n+1} \cap \{\omega \in \Omega|1\leq \|Y^N_n(\omega)\|\leq N^{1/2c}\}$ that
\begin{align}
	\label{eq_2_20}
	\|A_n\|^2 &\leq \frac{m^2}{4}\displaystyle\sum_{j_1,j_2=1}^{m}\|L^{j_1}g_{j_2}(Y^N_n)\|^2|\Delta W^{j_1}_n\Delta W^{j_2}_n - \delta_{j_1,j_2}h|^2  \notag\\
			&\leq \frac{m^2}{2}\displaystyle\sum_{j_1,j_2=1}^{m}(K + \|L^{j_1}g_{j_2}(0)\|)^2\|Y^N_n\|^2\cdot(|\Delta W^{j_1}_n\Delta W^{j_2}_n|^2 - \delta_{j_1,j_2}h^2)  \notag\\
			&\leq \frac{m^3}{2}\left(K + \max_{1\leq j_1,j_2\leq m}\|L^{j_1}g_{j_2}(0)\|\right)^2\left[\|\Delta W_n\|^2 + h^2\right]\|Y^N_n\|^2.
\end{align}
Inserting \eqref{eq_2_20} into \eqref{eq_2_19}, we obtain that on $\omega \in \Omega_{n+1} \cap \{\omega \in \Omega|1\leq \|Y^N_n(\omega)\|\leq N^{1/2c}\}$
\begin{align}
	\label{eq_2_21}
		\|Y^N_{n+1}\|^2 &\leq \|Y^N_n\|^2\bigg[1 + 4h^2(K + \|\phi(0)\|)^2 + 4h^2N(K + \|\varphi(0)\|)^2 + 4(K + \|g(0)\|)^2\|\Delta W_n\|^2   \notag\\
					  &\quad + 4h(K + \|\phi(0)\|)\bigg] + 2h(K + \|f(0)\|)\|Y^N_n\|^2 + 2m^3\left(K + \max_{1\leq j_1,j_2\leq m}\|L^{j_1}g_{j_2}(0)\|\right)^2  \notag\\
					  &\quad \times \left[\|\Delta W_n\|^2 + h^2\right]\|Y^N_n\|^2 + 2\langle Y^N_n, g(Y^N_n)\Delta W_n \rangle + 2\langle Y^N_n, A_n \rangle \notag\\
					  &\leq \|Y^N_n\|^2 \bigg [1 + \frac{2}{N}\bigg (2T^2(K + \|\phi(0)\|)^2 + 2T^2(K + \|\varphi(0)\|)^2 + 2T(K + \|\phi(0)\|)   \notag\\
					  &\quad + 2T(K + \|f(0)\|) + m^3T^2\left(K + \max_{1\leq j_1,j_2\leq m}\|L^{j_1}g_{j_2}(0)\|\right)^2\bigg )  \notag\\
					  &\quad + 2\left(2(K + \|g(0)\|)^2 + m^3\left(K + \max_{1\leq j_1,j_2\leq m}\|L^{j_1}g_{j_2}(0)\|\right)^2\right)\|\Delta W_n\|^2 + 2\alpha_k\bigg ]  \notag\\
					  &\leq \|Y^N_n\|^2\exp\left[\frac{2\lambda}{N} + 2\lambda\|\Delta W_n\|^2 + 2\alpha_k\right].
\end{align}
Then using the mathematical induction like the proof of Lemma 3.1 in \cite{hutzenthaler2012}, and combining \eqref{eq_2_10} and \eqref{eq_2_21}, we can establish \eqref{eq_2_9}.
\end{proof}

\begin{Lem}\label{lemma2_2}
	For $p \geq 1$, we have
	\begin{align}
		\label{eq_2_22}
		\sup_{N \in \mathbb{N}, N \geq 4\lambda pT} \mathbb{E}\left[\exp(p\lambda\sum_{k=0}^{N-1})\|\Delta W_k\|^2\right] < \infty.
	\end{align}
\end{Lem}

\begin{proof}
	The result is the same as Lemma 3.3 in \cite{hutzenthaler2012}, the only difference is the value of $\lambda$.
\end{proof}

\begin{Lem}\label{lemma2_3}
	 $\alpha_n$ is given by \eqref{eq_2_8}. Then for $p \geq 1$, we have
	\begin{align}
		\label{2_23}
		\sup_{z\in\{-1,1\}}\sup_{N\in\mathbb{N}}\mathbb{E}\left[\sup_{0\leq n\leq N}\exp\left(pz\displaystyle\sum_{k=0}^{n-1}\alpha_k\right)\right] < \infty.
	\end{align}
\end{Lem}

\begin{proof}
	The result is identical to Lemma 2.4 in \cite{wang2013}.
\end{proof}

\begin{Lem}\label{lemma2_4}
	Let $D_n$ be given by \eqref{eq_2_7}. Then for all $N \geq 8\lambda pT$ and $p \geq 1$
	\begin{align}
		\label{eq_2_24}
		\sup_{N \in\mathbb{N}}\mathbb{E}\left[\sup_{0\leq n\leq N}|D_n|^p\right] < \infty.
	\end{align}
\end{Lem}

\begin{proof}
$D_n$ here shares the same form as that in \cite{hutzenthaler2012}, but only with different $\lambda$ and $\alpha_n$. With the help of Lemma \ref{lemma2_2} and Lemma \ref{lemma2_3}, we can follow the proof of Lemma 3.5 in  \cite{hutzenthaler2012} to get the desired result.
\end{proof}

\begin{Lem}\label{lemma2_5}
	Let $\Omega_N$ be given by \eqref{eq_2_6} with $n = N.$ Then for all $p \geq 1$
	\begin{align}
		\label{eq_2_25}
		\sup_{N\in\mathbb{N}}(N^p\mathbb{P}[\Omega_{N}^c]) < \infty.
	\end{align}
\end{Lem}

\begin{proof}
	The proof is identical to the proof of Lemma 3.6 in \cite{hutzenthaler2012}.
\end{proof}

\begin{Lem}\label{lemma2_6}
	Let $k \in\mathbb{N}$ and $Z : [0,T]\times\Omega\rightarrow\mathbb{R}^{k\times m}$ be a predictable stochastic process satisfying $\mathbb{P}(\int_{0}^{T}\|Z_s\|^2ds < \infty) = 1.$ Then for all $t \in [0,T]$ and all $p \geq 2$
	\begin{align}
		\label{eq_2_26}
		\bigg \|\sup_{s\in[0,t]}\|\int_{0}^{s}Z_udW_u\|\bigg \|_{L^p(\Omega;\mathbb{R})} \leq p\left(\int_{0}^{t}\displaystyle\sum_{i=1}^{m}\|Z_se_i\|_{L^p(\Omega;\mathbb{R}^k)}^2ds\right)^{1/2}.
	\end{align}
Here, the vectors $e_i \in \mathbb{R}^m$ are orthogonal basis of the vector space $\mathbb{R}^m$.
\end{Lem}

\begin{proof}
	Combining Doob's maximal inequality and Lemma 7.7 in \cite{prato1992}, we can obtain the desired result.
\end{proof}

\begin{Lem}\label{lemma2_7}
	Let $k \in \mathbb{N}$ and let $Z_l : [0,T]\times\Omega\rightarrow\mathbb{R}^{k\times m}, l \in \{0,1,\cdots,N-1\}$ be a family of mappings such that $Z_l$ is $\mathcal{F}_{lT/N}/\mathcal{B}(\mathbb{R}^{k\times m})$-measurable. Then for $p\geq 2$ and $0\leq n\leq N$, we have
	\begin{align}
		\label{eq_2_27}
		\bigg \|\sup_{0\leq j\leq n}\|\displaystyle\sum_{l=0}^{j-1}Z_l\Delta W_l\|\bigg \|_{L^p(\Omega;\mathbb{R})} \leq p\left(\displaystyle\sum_{l=0}^{n-1}\displaystyle\sum_{i=1}^{m}\|Z_le_i\|_{L^p(\Omega;\mathbb{R}^k)}^2\frac{T}{N}\right)^{1/2}.
	\end{align}
\end{Lem}

This is a discrete version of the Burkholder-Davis-Gundy inequality \eqref{eq_2_26}.

\begin{Theo}\label{theorem2_1}
	$Y^N_n$ is given by \eqref{eq_1_4}. Then for $p \geq 1$, we have
	\begin{equation}
		\label{eq_2_28}
		\sup_{N \in\mathbb{N}}\left[\sup_{0\leq n\leq N}\mathbb{E}\|Y^N_n\|^p \right] < \infty.
	\end{equation}
\end{Theo}

\begin{proof}
	It follows from \eqref{eq_1_4} that
	\begin{align}
		\label{eq_2_29}
		\|Y^N_n\|_{L^p(\Omega;\mathbb{R}^d)} &\leq \|\xi\|_{L^p(\Omega;\mathbb{R}^d)} + \left\|\displaystyle\sum_{k=0}^{n-1}\phi(Y^N_k)h\right\|_{L^p(\Omega;\mathbb{R}^d)} + \left\|\displaystyle\sum_{k=0}^{n-1}\tilde{\varphi}(Y^N_k)h\right\|_{L^p(\Omega;\mathbb{R}^d)}  \notag\\
										   &\quad + \left\|\displaystyle\sum_{k=0}^{n-1}g(Y^N_k)\Delta W_k\right\|_{L^p(\Omega;\mathbb{R}^d)} + \left\|\displaystyle\sum_{k=0}^{n-1}A_k\right\|_{L^p(\Omega;\mathbb{R}^d)},
	\end{align}
where the notaion $A_n$ is given in \eqref{eq_2_12}. Employing triangle inequality and Lemma \ref{lemma2_7} we obtain
\begin{align}
	\label{eq_2_30}
	\left\|\displaystyle\sum_{k=0}^{n-1}g(Y^N_k)\Delta W_k\right\|_{L^p(\Omega;\mathbb{R}^d)} &\leq \left\|\displaystyle\sum_{k=0}^{n-1}(g(Y^N_k) - g(0))\Delta W_k\right\|_{L^p(\Omega;\mathbb{R}^d)} + \left\|\displaystyle\sum_{k=0}^{n-1}g(0)\Delta W_k\right\|_{L^p(\Omega;\mathbb{R}^d)}  \notag\\
																								 &\leq p\left(\frac{T}{N}\displaystyle\sum_{k=0}^{n-1}\displaystyle\sum_{i=0}^{m}\|g_i(Y^N_k) - g_i(0)\|_{L^p(\Omega;\mathbb{R}^d)}^2\right)^{1/2}  \notag\\
																								 &\quad + p\left(\frac{nT}{N}\displaystyle\sum_{i=0}^{m}\|g_i(0)\|_{L^p(\Omega;\mathbb{R}^d)}^2\right)^{1/2}  \notag\\
																								 &\leq p\left(mK^2h\displaystyle\sum_{k=0}^{n-1}\|Y^N_k\|_{L^p(\Omega;\mathbb{R}^d)}^2\right)^{1/2} + p\sqrt{Tm}\|g(0)\|.
\end{align}
From the independence of $Y^N_k$ and $\Delta W_k$, we obtain
\begin{align}
	\label{eq_2_31}
	\left\|\displaystyle\sum_{k=0}^{n-1}A_k\right\|_{L^p(\Omega;\mathbb{R}^d)} \leq &\frac{1}{2}\displaystyle\sum_{k=0}^{n-1}\displaystyle\sum_{j_1,j_2=1}^{m}\bigg[(K\|Y^N_k\|_{L^p(\Omega;\mathbb{R}^d)} + \|L^{j_1}g_{j_2}(0)\|) \notag\\
	&\times\|\Delta W^{j_1}_k\Delta W^{j_2}_k - \delta_{j_1,j_2}h\|_{L^p(\Omega;\mathbb{R}^d)}\bigg].
\end{align}
It follows from Lemma \ref{lemma2_7} and the mutual independence of $\Delta W_k$ that
\begin{equation}
	\label{eq_2_32}
	\frac{1}{2}\displaystyle\sum_{j_1,j_2=1}^{m}\|\Delta W^{j_1}_k\Delta W^{j_2}_k - \delta_{j_1,j_2}h\|_{L^p(\Omega;\mathbb{R}^d)} \leq c_{p,m}h,
\end{equation}
where $c_{p,m} = \frac{m^2 - m}{2}p^2 + 2mp^2 + \frac{m}{2}.$ Inserting \eqref{eq_2_32} into \eqref{eq_2_31} we obtain that
\begin{equation}
	\label{eq_2_33}
	\left\|\displaystyle\sum_{k=0}^{n-1}A_k\right\|_{L^p(\Omega;\mathbb{R}^d)} \leq Kc_{p,m}h\displaystyle\sum_{k=0}^{n-1}\|Y^N_k\|_{L^p(\Omega;\mathbb{R}^d)} + Tc_{p,m}\sup_{1\leq j_1,j_2\leq m}\|L^{j_1}g_{j_2}(0)\|.
\end{equation}
Additionally, we have
\begin{align}
	\label{eq_2_34}
	\left\|\displaystyle\sum_{k=0}^{n-1}\tilde{\varphi}(Y^N_k)h\right\|_{L^p(\Omega;\mathbb{R}^d)} &\leq N,  \\
	\label{eq_2_35}\left\|\displaystyle\sum_{k=0}^{n-1}\phi(Y^N_k)h\right\|_{L^p(\Omega;\mathbb{R}^d)}   	   &\leq hK\displaystyle\sum_{k=0}^{n-1}\|Y^N_k\|_{L^p(\Omega;\mathbb{R}^d)} + \frac{nT}{N}\|\phi(0)\|.
\end{align}
Combining \eqref{eq_2_29}-\eqref{eq_2_35} gives that
\begin{align}
	\label{eq_2_36}
	\|Y^N_n\|_{L^p(\Omega;\mathbb{R}^d)} &\leq \|\xi\|_{L^p(\Omega;\mathbb{R}^d)} + hK\displaystyle\sum_{k=0}^{n-1}\|Y^N_k\|_{L^p(\Omega;\mathbb{R}^d)} + \frac{nT}{N}\|\phi(0)\| + N  \notag\\
									   &\quad + Kc_{p,m}h\displaystyle\sum_{k=0}^{n-1}\|Y^N_k\|_{L^p(\Omega;\mathbb{R}^d)} + Tc_{p,m}\sup_{1\leq j_1,j_2\leq m}\|L^{j_1}g_{j_2}(0)\|  \notag\\
									   &\quad + p\left(mK^2h\displaystyle\sum_{k=0}^{n-1}\|Y^N_k\|_{L^p(\Omega;\mathbb{R}^d)}^2\right)^{1/2} + p\sqrt{Tm}\|g(0)\|.
\end{align}
Taking square on both sides gives that
\begin{align}
	\label{eq_2_37}
	\|Y^N_n\|^2_{L^p(\Omega;\mathbb{R}^d)} &\leq 4\left(\|\xi\|_{L^p(\Omega;\mathbb{R}^d)} + T\|\phi(0)\| + N + Tc_{p,m}\sup_{1\leq j_1,j_2\leq m}\|L^{j_1}g_{j_2}(0)\| +  p\sqrt{Tm}\|g(0)\|\right)^2 \notag\\
									   &\quad + 4K^2c_{p,m}^2Th\displaystyle\sum_{k=0}^{n-1}\|Y^N_k\|^2_{L^p(\Omega;\mathbb{R}^d)} + 4p^2mK^2h\displaystyle\sum_{k=0}^{n-1}\|Y^N_k\|_{L^p(\Omega;\mathbb{R}^d)}^2  \notag\\
									   &\quad + 4K^2Th\displaystyle\sum_{k=0}^{n-1}\|Y^N_k\|^2_{L^p(\Omega;\mathbb{R}^d)}.
\end{align}
Let $C_0 = T\|\phi(0)\| + Tc_{p,m}\displaystyle\sup_{1\leq j_1,j_2\leq m}\|L^{j_1}g_{j_2}(0)\| +  p\sqrt{Tm}\|g(0)\|$, and $C_1 = 2(K^2c_{p,m}^2T + p^2mK^2 + K^2T)T$, using Gronwall's lemma gives that
\begin{equation}
	\label{eq_2_38}
	\displaystyle\sup_{0\leq n\leq N}\|Y^N_n\|_{L^p(\Omega;\mathbb{R}^d)} \leq 2\exp(C_1)(\|\xi\|_{L^p(\Omega;\mathbb{R}^d)} + N + C_0).
\end{equation}
Employing H\"{o}lder's inequality, Lemma \ref{lemma2_5} and the estimate \eqref{eq_2_38}, we obtain that
\begin{align}
	\label{eq_2_39}
	\sup_{N \in\mathbb{N}}&\sup_{0\leq n\leq N}\|\mathds{1}_{(\Omega_{n})^c}Y^N_n\|_{L^p(\Omega;\mathbb{R}^d)}  \notag\\
	&\leq \sup_{N \in\mathbb{N}}\sup_{0\leq n\leq N}\bigg (\|\mathds{1}_{(\Omega_{n})^c}\|_{L^{2p}(\Omega;\mathbb{R}^d)}\|Y^N_n\|_{L^{2p}(\Omega;\mathbb{R}^d)}      \bigg )  \notag\\
    &\leq \bigg (\sup_{N \in\mathbb{N}}(N\cdot\|\mathds{1}_{(\Omega_{N})^c}\|_{L^{2p}(\Omega;\mathbb{R}^d)})\bigg )\cdot\bigg (\sup_{N \in\mathbb{N}}\sup_{0\leq n\leq N}\left(N^{-1}\|Y^N_n\|_{L^{2p}(\Omega;\mathbb{R}^d)}\right) \bigg )  \notag\\
    &\leq \bigg (\sup_{N \in\mathbb{N}}(N^{2p}\cdot\mathbb{P}[\Omega_{N}^c]\bigg )^{1/2p}\cdot\bigg (\sup_{N \in\mathbb{N}}\sup_{0\leq n\leq N}\left(N^{-1}\|Y^N_n\|_{L^{2p}(\Omega;\mathbb{R}^d)}\right) \bigg )  \notag\\
    &< \infty.
\end{align}
Furthermore, Lemma \ref{lemma2_1} and Lemma \ref{lemma2_4} imply that
\begin{equation}
	\label{eq_2_40}
	\sup_{N \in\mathbb{N}}\sup_{0\leq n\leq N}\|\mathds{1}_{\Omega_{n}}Y^N_n\|_{L^p(\Omega;\mathbb{R}^d)} \leq \sup_{N \in\mathbb{N}}\sup_{0\leq n\leq N}\|D_n\|_{L^p(\Omega;\mathbb{R}^d)} < \infty.
\end{equation}
Combining \eqref{eq_2_39} and \eqref{eq_2_40}, we can obtain the desired results.
\end{proof}

\section{Strong convergence order of the semi-tamed Milstein method}
\label{sec:StrongCon}
 In this section, the generic constant $C_{p,T}$ might vary in different places.

\begin{Assu}
	\label{assumption_3_1}
	Suppose that $f(x)$ and $g_i(x)$ are twice continuous differentiable. Furthermore, we suppose that there exist positive constants $K$, $q\geq 1$ such that for $x\in \mathbb{R}^d$ and $i \in \{1,2,\cdots,m\}$
	\begin{align}
		\label{eq_3_1}
		\|f''(x)\|_{L^{(2)}(\mathbb{R}^d;\mathbb{R}^d)} &\leq K(1 + \|x\|^q),  \\
		\label{eq_3_2}\|g_i''(x)\|_{L^{(2)}(\mathbb{R}^d;\mathbb{R}^d)} &\leq K,
	\end{align}
where, $\|\Phi''(x)\|_{L^{(2)}(\mathbb{R}^d;\mathbb{R}^d)} = \sup_{h_1,h_2 \in \mathbb{R},\|h_1\|\leq 1,\|h_2\|\leq 1}\|\Phi''(x)(h_1,h_2)\|$ for $\Phi:\mathbb{R}^d\rightarrow \mathbb{R}^d$, and $\Phi''(x):\mathbb{R}^d\times\mathbb{R}^d\rightarrow \mathbb{R}^d$ is a bilinear operator defined by \eqref{eq_3_8}.
\end{Assu}

In order to analyse the strong convergence of the semi-tamed Milstein method, we rewrite the scheme as
\begin{align}
	\label{Milstein_type}
	Y^N_{n+1} = Y^N_n& + h\phi(Y^N_n) + \frac{h\varphi(Y^N_n)}{1 + \|\varphi(Y^N_n)\|h} + g(Y^N_n)\Delta W_n \notag\\
		    & + \displaystyle\sum_{j_1,j_2=1}^{m}L^{j_1}g_{j_2}(Y^N_n)I_{j_1,j_2}^{t_n,t_{n+1}},
\end{align}
where $I_{j_1,j_2}^{t_n,t_{n+1}} = \int_{t_n}^{t_{n+1}}\int_{t_n}^{s2}dW_{s_1}^{j_1}dW_{s_2}^{j_2}$. Note that under commutativity condition \eqref{eq_1_5}, \eqref{eq_1_4} and \eqref{Milstein_type} are equivalent. Then we introduce time continuous interpolations of the time discrete numerical approximation, that is the time continuous approximation $\bar{Y}^N_s$ such that for $s\in[t_n, t_{n+1}]$
\begin{align}
	\label{eq_3_3}
	\bar{Y}^N_s &:= Y^N_n + (s - t_n)\phi(Y^N_n) + \frac{(s-t_n)\varphi(Y^N_n)}{1 + \|\varphi(Y^N_n)\|h} + g(Y^N_n)(W_s- W_{t_n}) + \displaystyle\sum_{j_1,j_2=1}^{m}L^{j_1}g_{j_2}(Y^N_n)I_{j_1,j_2}^{t_n,s}  \notag\\
			  &= Y^N_n + \int_{t_n}^{s}\phi(Y^N_n)du + \int_{t_n}^{s}\frac{\varphi(Y^N_n)}{1 + \|\varphi(Y^N_n)\|h}du + \displaystyle\sum_{i=1}^m\int_{t_n}^{s}g_i(Y^N_n)dW_u^i  \notag\\
			  &\quad + \displaystyle\sum_{i=1}^m\int_{t_n}^{s}\displaystyle\sum_{j=1}^mL^{j}g_{i}(Y^N_n)\Delta W_u^jdW_u^i  \notag\\
			  &= Y^N_n + \int_{t_n}^{s}\bar{f}(Y^N_n)du + \displaystyle\sum_{i=1}^m\int_{t_n}^{s}g_i(Y^N_n)dW_u^i + \displaystyle\sum_{i=1}^m\int_{t_n}^{s}\displaystyle\sum_{j=1}^mL^{j}g_{i}(Y^N_n)\Delta W_u^jdW_u^i,
\end{align}
where
\begin{align*}
	\bar{f}(x) &= f(x) - \frac{\|\varphi(x)\|\varphi(x)h}{1 + \|\varphi(x)\|h}, \\
	f(x) &= \phi(x) + \varphi(x),
\end{align*}
and
\begin{equation*}
	I_{j_1,j_2}^{t_n,s} = \int_{t_n}^{s}\int_{t_n}^{s2}dW_{s_1}^{j_1}dW_{s_2}^{j_2},\qquad\Delta W_s^j = \displaystyle\sum_{n=0}^\infty\mathds{1}_{\{t_n<s<t_{n+1}\}}\bigg (W_s^j - W_{t_n}^j\bigg ).
\end{equation*}
Obviously, $\bar{Y}^N_{t_n} = Y^N_n, n = 0, 1, \cdots, N$.
We can rewrite $\bar{Y}^N_t$ as
\begin{equation}
	\label{eq_3_4}
	\bar{Y}^N_t = Y^N_0 + \int_{0}^{t}\bar{f}(Y^N_{n_s})ds + \displaystyle\sum_{i=1}^m\int_{0}^{t}\left[g_i(Y^N_{n_s}) + \displaystyle\sum_{j=1}^mL^{j}g_{i}(Y^N_{n_s})\Delta W_s^j \right]dW_s^i,
\end{equation}
where
\begin{equation*}
	n_s := max\bigg\{n_s \in \{0, 1, 2,\cdots, N\}: t_{n_s} \leq s\bigg\},
\end{equation*}
which means the maximum integer such that $t_{n_s} \leq s$. Combining \eqref{eq_1_2} and \eqref{eq_3_4}, we obtain that
\begin{equation}
	\label{eq_3_5}
	X_t - \bar{Y}^N_t = \int_{0}^{t}[f(X_s) - \bar{f}(Y^N_{n_s})]ds + \displaystyle\sum_{i=1}^m\int_{0}^{t}\left[g_i(X_s) - g_i(Y^N_{n_s}) - \displaystyle\sum_{j=1}^mL^{j}g_{i}(Y^N_{n_s})\Delta W_s^j \right]dW_s^i.
\end{equation}

The Taylor's formula gives that: If a function $\Phi:\mathbb{R}^d\rightarrow \mathbb{R}^d$ is twice differentiable, then
\begin{equation}
	\label{eq_3_6}
	\Phi(\bar{Y}^N_s) - \Phi(Y^N_{n_s}) = \Phi'(Y^N_{n_s})(\bar{Y}^N_s - Y^N_{n_s}) + R_1(\Phi),
\end{equation}
where $R_1(\Phi)$ is the remainder term
\begin{equation}
	\label{eq_3_7}
	R_1(\Phi) = \int_{0}^{1}(1 - r)\Phi''(Y^N_{n_s} + r(\bar{Y}^N_s - Y^N_{n_s}))(\bar{Y}^N_s - Y^N_{n_s}, \bar{Y}^N_s - Y_{n_s})dr.
\end{equation}
Here we gives the derivative expression: for arbitrary $a, b_1, b_2 \in \mathbb{R}^d$
\begin{equation}
	\label{eq_3_8}
	\Phi'(a)(b_1) = \displaystyle\sum_{i=1}^{d}\frac{\partial \Phi}{\partial x^i}b_1^i, \quad \Phi''(a)(b_1, b_2) = \displaystyle\sum_{i,j=1}^{d}\frac{\partial^2 \Phi}{\partial x^i\partial x^j}b_1^ib_2^j.
\end{equation}
Replacing $\bar{Y}^N_s - Y^N_{n_s}$ in \eqref{eq_3_6} with \eqref{eq_3_3}, we obtain that
\begin{equation}
	\label{eq_3_9}
	\Phi(\bar{Y}^N_s) - \Phi(Y^N_{n_s}) = \Phi'(Y^N_{n_s})g(Y^N_{n_s})(W_s - W_{t_{n_s}}) + \bar{R}_1(\Phi),
\end{equation}
here
\begin{equation}
	\label{eq_3_10}
	\bar{R}_1(\Phi) = \Phi'(Y^N_{n_s})\left((s-t_{n_s})\bar{f}(Y^N_{n_s}) + \displaystyle\sum_{j_1,j_2=1}^mL^{j_1}g_{j_2}(Y^N_{n_s})I_{j_1,j_2}^{t_{n_s},s}\right) + R_1(\Phi).
\end{equation}
Replacing $\Phi$ by $g_i$ in \eqref{eq_3_9}, and considering $g_{i}'(x)g(x) = \displaystyle\sum_{j=1}^mL^{j}g_{i}(x)$, we obtain that
\begin{equation}
	\label{eq_3_11}
	\bar{R}_1(g_i) = g_i(\bar{Y}^N_s) - g_i(Y^N_{n_s}) - \displaystyle\sum_{j=1}^mL^{j}g_{i}(Y^N_{n_s})\Delta W_s^j.
\end{equation}

Before giving Theorem \ref{theorem3_1} which presents the convergence of the underlying numerical scheme, we first give the following five lemmas, which will be used to prove the theorem.
\begin{Lem}\label{lemma3_1}
	If Assumption \ref{assumption2_1} and the commutativity condition \eqref{eq_1_5} hold, then for $1\leq j_1,j_2\leq m$ and $p \geq 1$, we have the following estimates
	\begin{equation}
		\label{eq_3_13}
		\sup_{N \in\mathbb{N}}\sup_{0\leq n\leq N}\left[\mathbb{E}\|\phi(Y^N_n)\|^p \vee \mathbb{E}\|\varphi(Y^N_n)\|^p \vee \mathbb{E}\|f'(Y^N_n)\|^p \vee \mathbb{E}\|g(Y^N_n)\|^p \vee \mathbb{E}\|L^{j_1}g_{j_2}(Y^N_n)\|^p \right] < \infty.
	\end{equation}
\end{Lem}

\begin{proof}
	It can be obtained easily by Assumption \ref{assumption2_1} and Theorem \ref{theorem2_1}.
\end{proof}

\begin{Lem}\label{lemma3_2}
	Let Assumption \ref{assumption2_1} and Assumption \ref{assumption_3_1} hold, then for all $p \geq 1$
	\begin{equation}
		\label{eq_3_14}
		\sup_{t\in[0,T]}\left[\|X_t\|_{L^p(\Omega;\mathbb{R}^d)} \vee \|f(X_t)\|_{L^p(\Omega;\mathbb{R}^d)} \vee \|g(X_t)\|_{L^p(\Omega;\mathbb{R}^d)} \vee \|\bar{Y}^N_t\|_{L^p(\Omega;\mathbb{R}^d)}\right] < \infty.
	\end{equation}
\end{Lem}

\begin{proof}
	From conditions \eqref{eq_2_1} and \eqref{eq_2_2}, we obtain that the exact solution $X_t$ satisfies $\sup_{t\in[0,T]}\mathbb{E}\|X_t\|^p < \infty$ \cite{Mao1997}. The second estimate can be obtained by the globally Lipschitz condition on $\phi$, the polynomial condition on $\varphi'$ and $f(x)=\phi(x)+\varphi(x)$. The third estimate can be obtained by the globally Lipschitz condition on $g$. From \eqref{eq_3_3} and Lemma \ref{lemma3_1}, we can obtain the last estimate.
\end{proof}

\begin{Lem}\label{lemma3_3}
	Let Assumption \ref{assumption2_1}, Assumption \ref{assumption_3_1} and the commutativity condition \eqref{eq_1_5} hold. Then for $p \geq 1$, there exists a family of constants $C_{p,T} \geq 1$ such that
	\begin{equation}
		\label{eq_3_15}
		\sup_{t\in[0,T]}\|\bar{Y}^N_t - Y^N_{n_t}\|_{L^p(\Omega;\mathbb{R}^d)} \vee \sup_{t\in[0,T]}\|f(X_t) - f(X_{n_t})\|_{L^p(\Omega;\mathbb{R}^d)} \leq C_{p,T}h^{1/2}.
	\end{equation}
\end{Lem}

\begin{proof}
	Let $n_t$ be the maximum integer such that $t_{n_t} \leq t$. From \eqref{eq_3_3} we obtain
	\begin{equation}
		\label{eq_3_16}
		\bar{Y}^N_t - Y^N_{n_t} = (t - t_{n_t})\bar{f}(Y^N_{n_t}) + g(Y^N_{n_t})(W_t - W_{n_t}) + \displaystyle\sum_{j_1,j_2=1}^{m}L^{j_1}g_{j_2}I_{j_1,j_2}^{t_{n_t},t}.
	\end{equation}
The first estimate can be obtained follow a similar approach as the estimate of \eqref{eq_2_29}. One can obtain the second estimate similarly.
\end{proof}

\begin{Lem}\label{lemma3_4}
	Let $p \geq 1$ and let conditions in Assumption \ref{assumption2_1}, Assumption \ref{assumption_3_1} and the commutativity condition \eqref{eq_1_5} be fulfilled. Then for $i = \{1,2,\cdots,m\}$
	\begin{equation}
		\label{eq_3_17}
		\|\bar{R}(f)\|_{L^p(\Omega;\mathbb{R}^d)} \vee \|\bar{R}(g_i)\|_{L^p(\Omega;\mathbb{R}^d)} \leq C_{p,T}h.
	\end{equation}
\end{Lem}

\begin{proof}
	For the first estimate, considering equation \eqref{eq_3_9} and \eqref{eq_3_10}, and replacing $\Phi$ by $f$, we can obtain that
	\begin{equation}
		\label{eq_proof_3_4_1}
		\bar{R}(f) = f'(Y^N_{n_s})\left((s-t_{n_s})\bar{f}(Y^N_{n_s}) + \displaystyle\sum_{j_1,j_2=1}^mL^{j_1}g_{j_2}(Y^N_{n_s})I_{j_1,j_2}^{t_{n_s},s}\right) + R_1(f).
	\end{equation}
From \eqref{eq_3_7}, using H\"{o}lder inequality and Lemma 3.3, we can obtain the estimate of $R_1(f)$
\begin{align}
	\label{eq_proof_3_4_2}
	\|R_1(f)\|_{L^p(\Omega;\mathbb{R}^d)} &\leq \int_{0}^{1}(1 - r)\|f''(Y^N_{n_s} + r(\bar{Y}^N_s - Y^N_{n_s}))(\bar{Y}^N_s - Y^N_{n_s}, \bar{Y}^N_s - Y^N_{n_s})\|_{L^p(\Omega;\mathbb{R})}dr \notag\\
											&\leq \int_{0}^{1}\bigg\| \|f''(Y^N_{n_s} + r(\bar{Y}^N_s - Y^N_{n_s}))\|_{L^{(2)}(\mathbb{R}^d;\mathbb{R}^d)}\cdot\|\bar{Y}^N_s - Y^N_{n_s}\|^2\bigg\|_{L^p(\Omega;\mathbb{R})}dr \notag\\
											&\leq K\bigg\| (1 + \|Y^N_{n_s}\|^q + \|\bar{Y}^N_s\|^q)\cdot\|\bar{Y}^N_s - Y^N_{n_s}\|^2\bigg\|_{L^p(\Omega;\mathbb{R})} \notag\\
											&\leq K\|1 + \|Y^N_{n_s}\|^q + \|\bar{Y}^N_s\|^q\|_{L^{2p}(\Omega;\mathbb{R})}\cdot\|\bar{Y}^N_s - Y^N_{n_s}\|_{L^{4p}(\Omega;\mathbb{R})}^2  \notag\\
											&\leq C_{p,T}h.
\end{align}
The third inequality is derived by Jensen's inequality. Hence by \eqref{eq_proof_3_4_1} and the definition of $\bar{f}(x)$
\begin{align}
	\label{eq_proof_3_4_3}
	\|\bar{R}(f)\|_{L^p(\Omega;\mathbb{R}^d)} &\leq \|f'(Y^N_{n_s})(s-t_{n_s})\bar{f}(Y^N_{n_s})\|_{L^p(\Omega;\mathbb{R}^d)} + \|R_1(f)\|_{L^p(\Omega;\mathbb{R}^d)}   \notag\\
												&\quad + \bigg\| f'(Y^N_{n_s})\displaystyle\sum_{j_1,j_2=1}^mL^{j_1}g_{j_2}(Y^N_{n_s})I_{j_1,j_2}^{t_{n_s},s} \bigg\|_{L^p(\Omega;\mathbb{R}^d)} \notag\\
												&\leq h\bigg\|f'(Y^N_{n_s})\bigg( f(Y^N_{n_s}) - \frac{\|\varphi(Y^N_{n_s})\|\varphi(Y^N_{n_s})h}{1 + \|\varphi(Y^N_{n_s})\|h} \bigg)\bigg\|_{L^p(\Omega;\mathbb{R}^d)} \notag\\
												&\quad + \frac{1}{2}\displaystyle\sum_{j_1,j_2=1}^m\bigg\| f'(Y^N_{n_s})(L^{j_1}g_{j_2}(Y^N_{n_s})(\Delta W^{j_1}_s\Delta W^{j_2}_s - \delta_{j_1,j_2}h)) \bigg\|_{L^p(\Omega;\mathbb{R}^d)} \notag\\
												&\quad + \|R_1(f)\|_{L^p(\Omega;\mathbb{R}^d)}.
\end{align}
Using H\"{o}lder inequality, we derive from Lemma \ref{lemma3_1} that
\begin{align}
	\label{eq_proof_3_4_4}
	\bigg\|f'(Y^N_{n_s})&\bigg( f(Y^N_{n_s}) - \frac{\|\varphi(Y^N_{n_s})\|\varphi(Y^N_{n_s})h}{1 + \|\varphi(Y^N_{n_s})\|h} \bigg)\bigg\|_{L^p(\Omega;\mathbb{R}^d)}  \notag\\
	    				&\leq \|f'(Y^N_{n_s})f(Y^N_{n_s})\|_{L^p(\Omega;\mathbb{R}^d)} + \bigg\|f'(Y^N_{n_s})\frac{\|\varphi(Y^N_{n_s})\|\varphi(Y^N_{n_s})h}{1 + \|\varphi(Y^N_{n_s})\|h} \bigg)\bigg\|_{L^p(\Omega;\mathbb{R}^d)} \notag\\
	    				&\leq \|f'(Y^N_{n_s})\|_{L^{2p}(\Omega;\mathbb{R}^{d\times d})}\cdot\|f(Y^N_{n_s})\|_{L^{2p}(\Omega;\mathbb{R}^d)} + \|f'(Y^N_{n_s})\|_{L^{2p}(\Omega;\mathbb{R}^{d\times d})}\cdot\|\varphi(Y^N_{n_s})\|_{L^{2p}(\Omega;\mathbb{R}^d)} \notag\\
	    				& < \infty.
\end{align}

Similarly, by \eqref{eq_2_32} and Lemma 3.1, using H\"{o}lder inequality and the independence of $ Y^N_{n_s} $ and $ \Delta W^{j_1}_s, \Delta W^{j_2}_s $, there exists a constant $C_{p,T}$
\begin{align}
	\label{eq_proof_3_4_5}
	\frac{1}{2}\displaystyle\sum_{j_1,j_2=1}^m\bigg\| &f'(Y^N_{n_s})(L^{j_1}g_{j_2}(Y^N_{n_s})(\Delta W^{j_1}_s\Delta W^{j_2}_s - \delta_{j_1,j_2}h)) \bigg\|_{L^p(\Omega;\mathbb{R}^d)} \notag\\
													  &= \frac{1}{2}\displaystyle\sum_{j_1,j_2=1}^m\|f'(Y^N_{n_s})L^{j_1}g_{j_2}(Y^N_{n_s})\|_{L^p(\Omega;\mathbb{R}^d)}\times\|\Delta W^{j_1}_s\Delta W^{j_2}_s - \delta_{j_1,j_2}h\|_{L^p(\Omega;\mathbb{R}^d)}  \notag\\
													  &\leq C_{p,T}h.
\end{align}
Combining \eqref{eq_proof_3_4_3}, \eqref{eq_proof_3_4_4} and \eqref{eq_proof_3_4_5} we derive
\begin{equation}
	\|\bar{R}(f)\|_{L^p(\Omega;\mathbb{R}^d)} \leq C_{p,T}h.
\end{equation}
Following the similar approach, we can derive the second estimate
\begin{equation}
	\|\bar{R}(g_i)\|_{L^p(\Omega;\mathbb{R}^d)} \leq C_{p,T}h.
\end{equation}
\end{proof}

\begin{Lem}\label{lemma3_5}\cite{hutzenthaler2011b}
	Let $Z_1,\cdots,Z_N :\Omega\rightarrow\mathbb{R}$ be $\mathcal{F}/\mathcal{B}(\mathbb{R})$-measurable mappings with $\mathbb{E}\|Z_n\|^p < \infty$ for all $n \in \{1,\cdots,N\}$ and with $\mathbb{E}[Z_{n+1}|Z_1,\cdots,Z_n] = 0$ for all $n\in\{1,\cdots,N\}.$ Then
	\begin{equation}
		\label{eq_3_18}
		\|Z_1 + \cdots + Z_n\|_{L^p(\Omega;\mathbb{R})} \leq c_p\left(\|Z_1\|_{L^p(\Omega;\mathbb{R})}^2 + \cdots + \|Z_n\|_{L^p(\Omega;\mathbb{R})}^2\right)^{1/2},
	\end{equation}
for $p\geq 2$, where $c_p$ is independent of $n$, but dependent of $p$.
\end{Lem}

\begin{Theo}
	\label{theorem3_1}
	Suppose conditions in Assumption \ref{assumption2_1}, Assumption \ref{assumption_3_1} and the commutativity condition \eqref{eq_1_5} are fulfilled. Then, there exists a family of constants $C_{p,T} \geq 1$ and $p \in[1,\infty)$ such that
	\begin{equation}
		\label{eq_3_12}
		\left(\mathbb{E}\bigg[\sup_{t\in[0,T]}\|X_t - \bar{Y}^N_t\|^p \bigg] \right)^{1/p} \leq C_{p,T}\cdot h,\qquad h\in(0,1].
	\end{equation}
\end{Theo}

\begin{proof} Applying It\^{o}'s formula to \eqref{eq_3_5}, we obtain that
\begin{align}
	\label{eq_3_19}
	\|X_s - \bar{Y}^N_s\|^2 &= 2\int_{0}^{s}\langle X_u - \bar{Y}^N_u, f(X_u) - \bar{f}(Y^N_{n_u}) \rangle du  \notag\\
						  &\quad + 2\displaystyle\sum_{i=1}^m\int_{0}^{s}\langle X_u - \bar{Y}^N_{u}, g_i(X_u) - g_i(Y^N_{n_u}) - \displaystyle\sum_{j=1}^mL^{j}g_{i}(Y^N_{n_u})\Delta W_u^j \rangle dW_u^i  \notag\\
						  &\quad + \displaystyle\sum_{i=1}^m\int_{0}^{s}\bigg\|g_i(X_u) - g_i(Y^N_{n_u}) - \displaystyle\sum_{j=1}^mL^{j}g_{i}(Y^N_{n_u})\Delta W_u^j\bigg\|^2du  \notag\\
						  &= 2\int_{0}^{s}\langle X_u - \bar{Y}^N_u, f(X_u) - f(Y^N_{n_u}) \rangle du + 2\int_{0}^{s}\langle X_u - \bar{Y}^N_u, \frac{\|\varphi(Y^N_{n_u})\|\varphi(Y^N_{n_u})h}{1 + \|\varphi(Y^N_{n_u})\|h} \rangle du  \notag\\
						  &\quad + 2\displaystyle\sum_{i=1}^m\int_{0}^{s}\langle X_u - \bar{Y}^N_{u}, g_i(X_u) - g_i(Y^N_{n_u}) - \displaystyle\sum_{j=1}^mL^{j}g_{i}(Y^N_{n_u})\Delta W_u^j \rangle dW_u^i  \notag\\
						  &\quad + \displaystyle\sum_{i=1}^m\int_{0}^{s}\bigg\|g_i(X_u) - g_i(Y^N_{n_u}) - \displaystyle\sum_{j=1}^mL^{j}g_{i}(Y^N_{n_u})\Delta W_u^j\bigg\|^2du.
\end{align}
Employing the Cauchy-Schwarz inequality and the one-sided Lipschitz condition on $f$, we have that
\begin{align}
	\label{eq_3_20}
	\langle X_u - \bar{Y}^N_u, f(X_u) - f(Y^N_{n_u}) \rangle &\leq \langle X_u - \bar{Y}^N_u, f(X_u) - f(\bar{Y}^N_u) \rangle + \langle X_u - \bar{Y}^N_u, f(\bar{Y}^N_u) - f(Y^N_{n_u}) \rangle  \notag\\
															 &\leq K\|X_u - \bar{Y}^N_u\|^2 + \langle X_u - \bar{Y}^N_u, f(\bar{Y}^N_u) - f(Y^N_{n_u}) \rangle.
\end{align}
For the integrand of the second term in \eqref{eq_3_19}, we have
\begin{align}
	\label{eq_3_21}
	\langle X_u - \bar{Y}^N_u, \frac{\|\varphi(Y^N_{n_u})\|\varphi(Y^N_{n_u})h}{1 + \|\varphi(Y^N_{n_u})\|h} \rangle &\leq \frac{1}{2}\|X_u - \bar{Y}^N_u\|^2 + \frac{1}{2}\bigg \|\frac{\|\varphi(Y^N_{n_u})\|\varphi(Y^N_{n_u})h}{1 + \|\varphi(Y^N_{n_u})\|h}\bigg \|^2  \notag\\
																						   &\leq \frac{1}{2}\|X_u - \bar{Y}^N_u\|^2 + \frac{1}{2}h^2\|\varphi(Y^N_{n_u})\|^4.
\end{align}
Using \eqref{eq_2_2} and notation given by \eqref{eq_3_11}, we obtain
\begin{align}
	\label{eq_3_22}
	\bigg\|g_i(X_u) - g_i(Y^N_{n_u}) - \displaystyle\sum_{j=1}^mL^{j}g_{i}(Y^N_{n_u})\Delta W_u^j\bigg\|^2 &\leq 2K^2\|X_u - \bar{Y}^N_u\|^2 + 2\|\bar{R}_1(g_i)\|^2.
\end{align}
Inserting \eqref{eq_3_20}-\eqref{eq_3_22} into \eqref{eq_3_19}, we have
\begin{align}
	\label{eq_3_23}
	\|X_s - \bar{Y}^N_s\|^2 &\leq (2K + 1 + 2mK^2)\int_{0}^{s}\|X_u - \bar{Y}^N_u\|^2du + h^2\int_{0}^{s}\|\varphi(Y^N_{n_u})\|^4du  \notag\\
						  &\quad + 2\displaystyle\sum_{i=1}^m\int_{0}^{s}\|\bar{R}_1(g_i)\|^2du + 2\int_{0}^{s}\langle X_u - \bar{Y}^N_u, f(\bar{Y}^N_u) - f(Y^N_{n_u}) \rangle du  \notag\\
						  &\quad + 2\displaystyle\sum_{i=1}^m\int_{0}^{s}\langle X_u - \bar{Y}^N_{u}, g_i(X_u) - g_i(\bar{Y}^N_{u}) + \bar{R}_1(g_i) \rangle dW_u^i.
\end{align}
Therefore,
\begin{align}
	\label{eq_3_24}
	\sup_{0\leq s\leq t}\|X_s - \bar{Y}^N_s\|^2 &\leq (2K + 1 + 2mK^2)\int_{0}^{t}\|X_s - \bar{Y}^N_s\|^2ds + h^2\int_{0}^{t}\|\varphi(Y^N_{n_s})\|^4ds  \notag\\
											  &\quad + 2\displaystyle\sum_{i=1}^m\int_{0}^{t}\|\bar{R}_1(g_i)\|^2ds + 2\sup_{0\leq s\leq t}\int_{0}^{s}\langle X_u - \bar{Y}^N_u, f(\bar{Y}^N_u) - f(Y^N_{n_u}) \rangle du  \notag\\
											  &\quad + 2\sup_{0\leq s\leq t}\displaystyle\sum_{i=1}^m\int_{0}^{s}\langle X_u - \bar{Y}^N_{u}, g_i(X_u) - g_i(\bar{Y}^N_{u}) + \bar{R}_1(g_i) \rangle dW_u^i.
\end{align}
For $p \geq 4$, we have
\begin{align}
	\label{eq_3_25}
	\bigg\|\sup_{0\leq s\leq t}&\|X_s - \bar{Y}^N_s\|\bigg\|_{L^p(\Omega;\mathbb{R})}^2 = \bigg\|\sup_{0\leq s\leq t}\|X_s - \bar{Y}^N_s\|^2\bigg\|_{L^{p/2}(\Omega;\mathbb{R})}  \notag\\
																 &\leq (2K + 1 + 2mK^2)\int_{0}^{t}\|X_s - \bar{Y}^N_s\|_{L^p(\Omega;\mathbb{R}^d)}^2ds  \notag\\
																 &\quad + h^2\int_{0}^{t}\|\varphi(Y^N_{n_s})\|_{L^{2p}(\Omega;\mathbb{R}^d)}^4ds + 2\displaystyle\sum_{i=1}^m\int_{0}^{t}\|\bar{R}_1(g_i)\|_{L^{p}(\Omega;\mathbb{R}^d)}^2ds  \notag\\
																 &\quad + 2\bigg\|\sup_{0\leq s\leq t}\int_{0}^{s}\langle X_u - \bar{Y}^N_u, f(\bar{Y}^N_u) - f(Y^N_{n_u}) \rangle du\bigg\|_{L^{p/2}(\Omega;\mathbb{R})}  \notag\\
																 &\quad + 2\bigg\|\sup_{0\leq s\leq t}\displaystyle\sum_{i=1}^m\int_{0}^{s}\langle X_u - \bar{Y}^N_{u}, g_i(X_u) - g_i(\bar{Y}^N_{u}) + \bar{R}_1(g_i) \rangle dW_u^i\bigg\|_{L^{p/2}(\Omega;\mathbb{R})}.
\end{align}
It's easy to obtain that
\begin{align}
	\label{eq_3_26}
	2\bigg\|&\sup_{0\leq s\leq t}\displaystyle\sum_{i=1}^m\int_{0}^{s}\langle X_u - \bar{Y}^N_{u}, g_i(X_u) - g_i(\bar{Y}^N_{u}) + \bar{R}_1(g_i) \rangle dW_u^i\bigg\|_{L^{p/2}(\Omega;\mathbb{R})} \notag\\
		    &\leq \frac{1}{4}\sup_{0\leq s\leq t}\|X_s - \bar{Y}^N_s\|_{L^{p}(\Omega;\mathbb{R}^d)}^2 + 2p^2m^2K\int_{0}^{t}\|X_s - \bar{Y}^N_s\|_{L^p(\Omega;\mathbb{R}^d)}^2ds  \notag\\
		    &\quad + 2p^2m\displaystyle\sum_{i=1}^m\int_{0}^{t}\|\bar{R}_1(g_i)\|_{L^{p}(\Omega;\mathbb{R}^d)}^2ds,
\end{align}
and
\begin{align}
	\label{eq_3_27}
	2\bigg\|&\sup_{0\leq s\leq t}\int_{0}^{s}\langle X_u - \bar{Y}^N_u, f(\bar{Y}^N_u) - f(Y^N_{n_u}) \rangle du\bigg\|_{L^{p/2}(\Omega;\mathbb{R})}  \notag\\
			&= 2\bigg\|\sup_{0\leq s\leq t}\int_{0}^{s}\langle X_u - \bar{Y}^N_u, f'(Y^N_{n_u})(g(Y^N_{n_u})\Delta W_u) + \bar{R}_1(f) \rangle du\bigg\|_{L^{p/2}(\Omega;\mathbb{R})}  \notag\\
			&\leq M + 2\bigg\|\sup_{0\leq s\leq t}\int_{0}^{s}\langle X_u - \bar{Y}^N_u, \bar{R}_1(f) \rangle du\bigg\|_{L^{p/2}(\Omega;\mathbb{R})}  \notag\\
			&\leq M + \int_{0}^{t}\|X_s - \bar{Y}^N_s\|_{L^{p}(\Omega;\mathbb{R}^d)}^2ds + \int_{0}^{t}\|\bar{R}_1(f)\|_{L^{p}(\Omega;\mathbb{R}^d)}^2ds,
\end{align}
where
\begin{equation}
	\label{eq_3_28}
	M = 2\bigg\|\sup_{0\leq s\leq t}\int_{0}^{s}\langle X_u - \bar{Y}^N_u, f'(Y^N_{n_u})g(Y^N_{n_u})\Delta W_u \rangle du\bigg\|_{L^{p/2}(\Omega;\mathbb{R})}.
\end{equation}
Furthermore, we have
\begin{align}
	\label{eq_3_29}
	X_u - \bar{Y}^N_u &= X_{t_{n_u}} - Y^N_{n_u} + \displaystyle\int_{t_{n_u}}^{u}f(X_r)dr - \displaystyle\int_{t_{n_u}}^{u}\bar{f}(Y^N_{n_u})dr  \notag\\
			  &\quad + \displaystyle\sum_{i=1}^{m}\displaystyle\int_{t_{n_u}}^{u}\left[ g_{i}(X_r) - g_{i}(Y^N_{n_r}) - \displaystyle\sum_{j=1}^{m}L^jg_{i}(Y^N_{n_r})\Delta W_{r}^j \right]dW_{r}^i	\notag\\
			  &= \displaystyle\int_{t_{n_u}}^{u}[f(X_r) - f(X_{t_{n_u}})]dr +  \displaystyle\sum_{i=1}^{m}\displaystyle\int_{t_{n_u}}^{u}[g_{i}(X_r) - g_{i}(\bar{Y}^N_r)]dW_{r}^i   \notag\\
			  &\quad + \displaystyle\sum_{i=1}^{m}\displaystyle\int_{t_{n_u}}^{u}\bar{R}_1(g_{i})dW_{r}^i + (u - t_{n_u})f(X_{t_{n_u}}) - (u - t_{n_u})\bar{f}(Y^N_{n_u}) + X_{t_{n_u}} - Y^N_{n_u}.
\end{align}
Replacing $X_u - \bar{Y}^N_u$ in \eqref{eq_3_28} by \eqref{eq_3_29}, we obtain
\begin{align}
	\label{eq_3_30}
	M &\leq 2\bigg\|\sup_{0\leq s\leq t}\int_{0}^{s}\left\langle \int_{t_{n_u}}^{u}[f(X_r) - f(X_{t_{n_u}})]dr, f'(Y^N_{n_u})g(Y^N_{n_u})\Delta W_u \right\rangle du\bigg\|_{L^{p/2}(\Omega;\mathbb{R})}  \notag\\
	  &\quad + 2\bigg\|\sup_{0\leq s\leq t}\int_{0}^{s}\left\langle \displaystyle\sum_{i=1}^{m}\displaystyle\int_{t_{n_u}}^{u}[g_{i}(X_r) - g_{i}(\bar{Y}^N_r)]dW_{r}^i, f'(Y^N_{n_u})g(Y^N_{n_u})\Delta W_u \right\rangle du\bigg\|_{L^{p/2}(\Omega;\mathbb{R})}  \notag\\
	  &\quad + 2\bigg\|\sup_{0\leq s\leq t}\int_{0}^{s}\left\langle \displaystyle\sum_{i=1}^{m}\displaystyle\int_{t_{n_u}}^{u}\bar{R}_1(g_{i})dW_{r}^i, f'(Y^N_{n_u})g(Y^N_{n_u})\Delta W_u \right\rangle du\bigg\|_{L^{p/2}(\Omega;\mathbb{R})}  \notag\\
	  &\quad + 2\bigg\|\sup_{0\leq s\leq t}\int_{0}^{s}\left\langle \zeta_{n_u}, f'(Y^N_{n_u})g(Y^N_{n_u})\Delta W_u \right\rangle du\bigg\|_{L^{p/2}(\Omega;\mathbb{R})}  \notag\\
	  &\quad + 2\bigg\|\sup_{0\leq s\leq t}\int_{0}^{s}\left\langle X_{t_{n_u}} - Y^N_{n_u}, f'(Y^N_{n_u})g(Y^N_{n_u})\Delta W_u \right\rangle du\bigg\|_{L^{p/2}(\Omega;\mathbb{R})}  \notag\\
	  &:= M_1 + M_2 + M_3 +M_4 + M_5,
\end{align}
where $\zeta_{n_u} \in \mathcal{F}_{t_{n_u}}$ is defined by
\begin{equation}
	\label{eq_3_31}
	\zeta_{n_u} = (u - t_{n_u})f(X_{t_{n_u}}) - (u - t_{n_u})\bar{f}(Y^N_{n_u}).
\end{equation}
By H\"{o}lder inequality, Lemma \ref{lemma3_1} and Lemma \ref{lemma3_5}, we can establish the following estimate
\begin{equation}
	\label{eq_3_32}
	\|f'(Y^N_{n_u})g(Y^N_{n_u})\Delta W_u\|_{L^{p}(\Omega;\mathbb{R}^d)} \leq C_{p,T}h^{1/2}.
\end{equation}
By Lemma \ref{lemma3_3} and \eqref{eq_3_32}, we have
\begin{align}
	\label{eq_3_33}
	M_1 &\leq 2\int_{0}^{t}\int_{t_{n_u}}^{u}\|f(X_r) - f(X_{t_{n_u}})\|_{L^{p}(\Omega;\mathbb{R}^d)}\cdot\|f'(Y^N_{n_u})g(Y^N_{n_u})\Delta W_u\|_{L^{p}(\Omega;\mathbb{R}^d)}drdu \notag\\
	    &\leq C_{p,T}h^2.
\end{align}
By H\"{o}lder inequality, \eqref{eq_2_2} and \eqref{eq_2_26}, we obtain
\begin{align}
	\label{eq_3_34}
	M_2 &\leq 2\int_{0}^{t}\bigg\| \displaystyle\sum_{i=1}^{m}\displaystyle\int_{t_{n_u}}^{u}[g_{i}(X_r) - g_{i}(\bar{Y}^N_r)]dW_{r}^i \bigg\|_{L^{p}(\Omega;\mathbb{R}^d)}\cdot\|f'(Y^N_{n_u})g(Y^N_{n_u})\Delta W_u\|_{L^{p}(\Omega;\mathbb{R}^d)}du  \notag\\
	    &\leq \int_{0}^{t}\frac{1}{h}\bigg\| \displaystyle\sum_{i=1}^{m}\displaystyle\int_{t_{n_u}}^{u}[g_{i}(X_r) - g_{i}(\bar{Y}^N_r)]dW_{r}^i \bigg\|_{L^{p}(\Omega;\mathbb{R}^d)}^2du  \notag\\
	    &\quad + \int_{0}^{t}h\|f'(Y^N_{n_u})g(Y^N_{n_u})\Delta W_u\|_{L^{p}(\Omega;\mathbb{R}^d)}^2du  \notag\\
	    &\leq \frac{p^2}{h}\int_{0}^{t}\int_{t_{n_u}}^{u}\displaystyle\sum_{i=1}^{m}\|g_{i}(X_r) - g_{i}(\bar{Y}^N_r)\|_{L^{p}(\Omega;\mathbb{R}^d)}^2drdu + C_{p,T}h^2  \notag\\
	    &\leq mp^2K^2\int_{0}^{t}\sup_{0\leq r\leq u}\|X_r - \bar{Y}^N_r\|_{L^{p}(\Omega;\mathbb{R}^d)}^2du + C_{p,T}h^2.
\end{align}
Similarly, we have the following estimate for $M_3$
\begin{align}
	\label{eq_3_35}
	M_3 &\leq 2\int_{0}^{t}\bigg\| \displaystyle\sum_{i=1}^{m}\displaystyle\int_{t_{n_u}}^{u}\bar{R}_1(g_{i})dW_{r}^i \bigg\|_{L^{p}(\Omega;\mathbb{R}^d)}\cdot\|f'(Y^N_{n_u})g(Y^N_{n_u})\Delta W_u\|_{L^{p}(\Omega;\mathbb{R}^d)}du  \notag\\
	    &\leq 2p\int_{0}^{t}\left( \displaystyle\int_{t_{n_u}}^{u}\displaystyle\sum_{i=1}^{m}\|\bar{R}_1(g_{i})\|_{L^{p}(\Omega;\mathbb{R}^d)}^2dr \right)^{1/2}\cdot\|f'(Y^N_{n_u})g(Y^N_{n_u})\Delta W_u\|_{L^{p}(\Omega;\mathbb{R}^d)}du  \notag\\
	    &\leq C_{p,T}h^2.
\end{align}
Furthermore,
\begin{align}
	\label{eq_3_36}
	M_4 &\leq 2\bigg\|\sup_{0\leq s\leq t} \bigg|\displaystyle\sum_{k=0}^{n_s-1}\int_{t_k}^{t_{k+1}}\left\langle \zeta_k, f'(Y^N_k)g(Y^N_k)\Delta W_u \right\rangle du\bigg| \bigg\|_{L^{p/2}(\Omega;\mathbb{R})} \notag\\
	    &\quad + 2\bigg\|\sup_{0\leq s\leq t} \bigg|\int_{t_{n_s}}^{s}\left\langle \zeta_{n_s}, f'(Y^N_{n_s})g(Y^N_{n_s})\Delta W_u \right\rangle du\bigg| \bigg\|_{L^{p/2}(\Omega;\mathbb{R})}  \notag\\
	    &:= M_{41} + M_{42}.
\end{align}
Obviously, the discrete time process
\begin{equation}
	\label{eq_3_37}
	\chi_n := \left\{ \displaystyle\sum_{k=0}^{n-1}\int_{t_k}^{t_{k+1}}\left\langle \zeta_k, f'(Y^N_k)g(Y^N_k)\Delta W_u \right\rangle du \right\},\quad n\in\{0,1,\cdots,N\}
\end{equation}
is an $\{\mathcal{F}_{t_n} : 0\leq n\leq N\}$-martingale.
Hence by Doob's maximal inequality, H\"{o}lder inequality and Lemma \ref{lemma3_5}, we obtain that for $p\geq 4$
\begin{align}
	\label{eq_3_38}
	M_{41} &\leq \frac{2p}{p - 2}\bigg\| \displaystyle\sum_{k=0}^{n_t-1}\int_{t_k}^{t_{k+1}}\left\langle \zeta_k, f'(Y^N_k)g(Y^N_k)\Delta W_u \right\rangle du \bigg\|_{L^{p/2}(\Omega;\mathbb{R})}  \notag\\
	       &\leq \frac{2pc_{p/2}}{p - 2}\left( \displaystyle\sum_{k=0}^{n_t-1}\bigg\| \int_{t_k}^{t_{k+1}}\left\langle \zeta_k, f'(Y^N_k)g(Y^N_k)\Delta W_u \right\rangle du \bigg\|_{L^{p/2}(\Omega;\mathbb{R})}^2 \right)^{1/2}  \notag\\
	       &\leq \frac{2pc_{p/2}}{p - 2}\left( \displaystyle\sum_{k=0}^{n_t-1}h \int_{t_k}^{t_{k+1}}\bigg\|\left\langle \zeta_k, f'(Y^N_k)g(Y^N_k)\Delta W_u \right\rangle \bigg\|_{L^{p/2}(\Omega;\mathbb{R})}^2du \right)^{1/2}  \notag\\
	       &\leq \frac{2pc_{p/2}}{p - 2}\left( \displaystyle\sum_{k=0}^{n_t-1}h \int_{t_k}^{t_{k+1}} \|\zeta_k\|_{L^{p}(\Omega;\mathbb{R}^d)}^2\cdot\|f'(Y^N_k)g(Y^N_k)\Delta W_u\|_{L^{p}(\Omega;\mathbb{R}^d)}^2du \right)^{1/2}.
\end{align}
From Lemma \ref{lemma3_1} and Lemma \ref{lemma3_2}, we have
\begin{equation}
	\label{eq_3_39}
	\|\zeta_k\|_{L^{p}(\Omega;\mathbb{R}^d)} \leq h\|f(X_{t_k})\|_{L^{p}(\Omega;\mathbb{R}^d)} + h\|\bar{f}(Y^N_k)\|_{L^{p}(\Omega;\mathbb{R}^d)}  \leq C_{p,T}h.
\end{equation}
Combining \eqref{eq_3_32}, \eqref{eq_3_38} and \eqref{eq_3_39}, we have
\begin{equation}
	\label{eq_3_40}
	M_{41} \leq C_{p,T}h^2.
\end{equation}
Furthermore, by H\"{o}lder inequality, \eqref{eq_3_32} and \eqref{eq_3_39}, we obtain that for $p\geq 4$
\begin{align}
	\label{eq_3_41}
	\left( \frac{M_{42}}{2} \right)^{p/2} &= \mathbb{E}\left( \sup_{0\leq s\leq t} \bigg|\int_{t_{n_s}}^{s}\left\langle \zeta_{n_s}, f'(Y^N_{n_s})g(Y^N_{n_s})\Delta W_u \right\rangle du\bigg| \right)^{p/2}  \notag\\
										  &\leq h^{p/2 - 1}\mathbb{E}\left( \sup_{0\leq s\leq t}\int_{t_{n_s}}^{s}\bigg|\left\langle \zeta_{n_s}, f'(Y^N_{n_s})g(Y^N_{n_s})\Delta W_u \right\rangle\bigg|^{p/2}du \right)  \notag\\
										  &\leq h^{p/2 - 1}\mathbb{E}\bigg( \displaystyle\sum_{k=0}^{n_t-1}\int_{t_k}^{t_{k+1}}\bigg|\left\langle \zeta_{k}, f'(Y^N_{k})g(Y^N_{k})\Delta W_u \right\rangle\bigg|^{p/2}du   \notag\\
										  &\quad + \int_{t_{n_t}}^{t}\bigg|\left\langle \zeta_{n_t}, f'(Y^N_{n_t})g(Y^N_{n_t})\Delta W_u \right\rangle\bigg|^{p/2}du\bigg)  \notag\\
										  &= h^{p/2 - 1}\int_{0}^{t}\mathbb{E}|\left\langle \zeta_{n_u}, f'(Y^N_{n_u})g(Y^N_{n_u})\Delta W_u \right\rangle|^{p/2}du  \notag\\
										  &\leq h^{p/2 - 1}\int_{0}^{t}\|\zeta_{n_u}\|_{L^{p}(\Omega;\mathbb{R})}^{p/2}\cdot\|f'(Y^N_{n_u})g(Y^N_{n_u})\Delta W_u\|_{L^{p}(\Omega;\mathbb{R}^d)}^{p/2}du  \notag\\
										  &\leq C_{p,T}h^{\frac{5p}{4} - 1}.
\end{align}
From \eqref{eq_3_41}, we obtain that, there exists a suitable constant $C_{p,T}$ such that
\begin{equation}
	\label{eq_3_42}
	M_{42} \leq 2C_{p.T}^{2/p}h^{(5p-4)/2p} \leq C_{p,T}h^2,
\end{equation}
for $p\geq 4$ and $h\in ( 0, 1 ]$.
Combining \eqref{eq_3_36}, \eqref{eq_3_40} and \eqref{eq_3_42}, we have that for $p\geq 4$
\begin{equation}
	\label{eq_3_43}
	M_4 \leq C_{p,T}h^2.
\end{equation}
Similarly, we split $M_5$ into two terms
\begin{align}
	\label{eq_3_44}
	M_5 &\leq 2\bigg\|\sup_{0\leq s\leq t} \bigg|\displaystyle\sum_{k=0}^{n_s-1}\int_{t_k}^{t_{k+1}}\left\langle X_{t_k} - Y^N_k, f'(Y^N_k)g(Y^N_k)\Delta W_u \right\rangle du\bigg| \bigg\|_{L^{p/2}(\Omega;\mathbb{R})} \notag\\
		&\quad + 2\bigg\|\sup_{0\leq s\leq t} \bigg|\int_{t_{n_s}}^{s}\left\langle X_{t_s} - Y^N_{t_{n_s}}, f'(Y^N_{n_s})g(Y^N_{n_s})\Delta W_u \right\rangle du\bigg| \bigg\|_{L^{p/2}(\Omega;\mathbb{R})}  \notag\\
		&:= M_{51} + M_{52}.
\end{align}
Using the fact that $\bar{Y}^N_{t_k} = Y^N_k, k = 0,1,\cdots,N-1$, we obtain the following estimate for $M_{51}$
\begin{align}
	\label{eq_3_45}
	M_{51} &\leq \frac{2pc_{p/2}}{p - 2}\left( \displaystyle\sum_{k=0}^{n_t-1}h \int_{t_k}^{t_{k+1}} \|X_{t_k} - Y^N_k\|_{L^{p}(\Omega;\mathbb{R}^d)}^2\cdot\|f'(Y^N_k)g(Y^N_k)\Delta W_u\|_{L^{p}(\Omega;\mathbb{R}^d)}^2du \right)^{1/2}  \notag\\
		   &\leq \sup_{0\leq s\leq t}\|X_s - \bar{Y}^N_s\|_{L^{p}(\Omega;\mathbb{R}^d)}\cdot\frac{2pc_{p/2}}{p - 2}\left( \displaystyle\sum_{k=0}^{n_t-1}h \int_{t_k}^{t_{k+1}}\|f'(Y^N_k)g(Y^N_k)\Delta W_u\|_{L^{p}(\Omega;\mathbb{R}^d)}^2du \right)^{1/2}  \notag\\
		   &\leq \frac{1}{4}\sup_{0\leq s\leq t}\|X_s - \bar{Y}^N_s\|_{L^{p}(\Omega;\mathbb{R}^d)}^2 + \frac{4p^2c_{p/2}^2}{(p - 2)^2}\displaystyle\sum_{k=0}^{n_t-1}h \int_{t_k}^{t_{k+1}}\|f'(Y^N_k)g(Y^N_k)\Delta W_u\|_{L^{p}(\Omega;\mathbb{R}^d)}^2du  \notag\\
		   &\leq \frac{1}{4}\sup_{0\leq s\leq t}\|X_s - \bar{Y}^N_s\|_{L^{p}(\Omega;\mathbb{R}^d)}^2 + C_{p,T}h^2.
\end{align}
Furthermore,
\begin{align}
	\label{eq_3_46}
	\left( \frac{M_{52}}{2} \right)^{p/2} &\leq h^{p/2 - 1}\int_{0}^{t}\|X_{t_{n_u}} - Y^N_{n_u}\|_{L^{p}(\Omega;\mathbb{R}^d)}^{p/2}\cdot\|f'(Y^N_{n_u})g(Y^N_{n_u})\Delta W_u\|_{L^{p}(\Omega;\mathbb{R}^d)}^{p/2}du  \notag\\
										  &\leq \sup_{0\leq s\leq t}\|X_s - \bar{Y}^N_s\|_{L^{p}(\Omega;\mathbb{R}^d)}^{p/2}\cdot C_{p,T}h^{(3p/4) - 1}.
\end{align}
Therefore
\begin{align}
	\label{eq_3_47}
	M_{52} &\leq 2\sup_{0\leq s\leq t}\|X_s - \bar{Y}^N_s\|_{L^{p}(\Omega;\mathbb{R}^d)}\cdot C_{p,T}^{2/p}h^{(3p-4)/2p}  \notag\\
		   &\leq \frac{1}{4}\sup_{0\leq s\leq t}\|X_s - \bar{Y}^N_s\|_{L^{p}(\Omega;\mathbb{R}^d)}^2 + 4C_{p,T}^{4/p}h^{(3p-4)/p} \notag\\
		   &\leq \frac{1}{4}\sup_{0\leq s\leq t}\|X_s - \bar{Y}^N_s\|_{L^{p}(\Omega;\mathbb{R}^d)}^2 + 4C_{p,T}^{4/p}h^{(3p-4)/p},
\end{align}
which implies that for $p\geq 4$ and $h\in ( 0, 1 ]$ there exists a suitable constant $C_{p,T}$ such that
\begin{equation}
	\label{eq_3_48}
	M_{52} \leq \frac{1}{4}\sup_{0\leq s\leq t}\|X_s - \bar{Y}^N_s\|_{L^{p}(\Omega;\mathbb{R}^d)}^2 + C_{p,T}h^2.
\end{equation}
Combining \eqref{eq_3_44}, \eqref{eq_3_45} and \eqref{eq_3_48} we derive that
\begin{equation}
	\label{eq_3_49}
	M_5 \leq \frac{1}{2}\sup_{0\leq s\leq t}\|X_s - \bar{Y}^N_s\|_{L^{p}(\Omega;\mathbb{R}^d)}^2 + C_{p,T}h^2.
\end{equation}
Inserting \eqref{eq_3_33}, \eqref{eq_3_34}, \eqref{eq_3_35}, \eqref{eq_3_43} and \eqref{eq_3_49} into \eqref{eq_3_30} yields
\begin{equation}
	\label{eq_3_50}
	M \leq mp^2K^2\int_{0}^{t}\sup_{0\leq r\leq u}\|X_r - \bar{Y}^N_r\|_{L^{p}(\Omega;\mathbb{R}^d)}^2du + \frac{1}{2}\sup_{0\leq s\leq t}\|X_s - \bar{Y}^N_s\|_{L^{p}(\Omega;\mathbb{R}^d)}^2 + C_{p,T}h^2.
\end{equation}
Hence, combining \eqref{eq_3_24}-\eqref{eq_3_28} we obtain that
\begin{align}
	\label{eq_3_51}
	\frac{3}{4}\bigg\|\sup_{0\leq s\leq t}\|X_s - \bar{Y}^N_s\|\bigg\|_{L^p(\Omega;\mathbb{R})}^2 &\leq 2(K + 1 + mK^2 + p^2m^2K + \frac{1}{2}mp^2K^2)  \notag\\
	 																				 &\quad\times\int_{0}^{t}\sup_{0\leq u\leq s}\|X_u - \bar{Y}^N_u\|_{L^{p}(\Omega;\mathbb{R}^d)}^2ds \notag\\
																					 &\quad + \frac{1}{2}\sup_{0\leq s\leq t}\|X_s - \bar{Y}^N_s\|_{L^{p}(\Omega;\mathbb{R}^d)}^2 + C_{p,T}h^2  \notag\\
																					 &\quad + 2(1 + p^2m)\displaystyle\sum_{i=1}^m\int_{0}^{t}\|\bar{R}_1(g_i)\|_{L^{p}(\Omega;\mathbb{R}^d)}^2ds \notag\\
																					 &\quad + h^2\int_{0}^{t}\|\varphi(Y^N_{n_s})\|_{L^{2p}(\Omega;\mathbb{R}^d)}^4ds + \int_{0}^{t}\|\bar{R}_1(f)\|_{L^{p}(\Omega;\mathbb{R}^d)}^2ds.																					
\end{align}
Additionally, from Lemma \ref{lemma3_1} and Lemma \ref{lemma3_4}, we can obtain that
\begin{align}
	\label{eq_3_52}
	2(1 + p^2m)\displaystyle\sum_{i=1}^m\int_{0}^{t}\|\bar{R}_1(g_i)\|_{L^{p}(\Omega;\mathbb{R}^d)}^2ds &\leq 2(1 + p^2m)\int_{0}^{t}C_{p,T}^2h^2ds \leq C_{p,T}h^2, \\
											 \label{eq_3_53}h^2\int_{0}^{t}\|\varphi(Y^N_{n_s})\|_{L^{2p}(\Omega;\mathbb{R}^d)}^4ds &\leq C_{p,T}h^2,  \\
											 \label{eq_3_54}\int_{0}^{t}\|\bar{R}_1(f)\|_{L^{p}(\Omega;\mathbb{R}^d)}^2ds &\leq C_{p,T}h^2.
\end{align}
Thus, we can get by \eqref{eq_3_51} that
\begin{equation}
	\label{eq_3_55}
	\bigg\|\sup_{0\leq s\leq t}\|X_s - \bar{Y}^N_s\|\bigg\|_{L^p(\Omega;\mathbb{R})}^2 \leq C_{p,T}\int_{0}^{t}\bigg\|  \sup_{0\leq u\leq s}\|X_u - \bar{Y}^N_u\|\bigg\|_{L^{p}(\Omega;\mathbb{R})}^2ds + C_{p,T}h^2.
\end{equation}
The Gronwall inequality gives the desired result for $p\geq 4.$ The assertion for $1\leq p<4$ can be obtained by using H\"{o}lder inequality. To be precisely, for $1\leq p<4$, using H\"{o}lder inequality, we have
\begin{align}
	\bigg\|\sup_{0\leq s\leq t}\|X_s - \bar{Y}^N_s\|\bigg\|_{L^p(\Omega;\mathbb{R})}^2 &= \bigg\|\sup_{0\leq s\leq t}\|X_s - \bar{Y}^N_s\|^2\bigg\|_{L^{p/2}(\Omega;\mathbb{R})} \notag\\
																					  &\leq \|1\|_{L^{p}(\Omega;\mathbb{R})}\times \bigg\|\sup_{0\leq s\leq t}\|X_s - \bar{Y}^N_s\|^2\bigg\|_{L^{p}(\Omega;\mathbb{R})} \notag\\
																					  & \leq \cdots \notag\\
																					  &\leq \bigg\|\sup_{0\leq s\leq t}\|X_s - \bar{Y}^N_s\|^2\bigg\|_{L^{2np}(\Omega;\mathbb{R})}.
\end{align}
This process can be continued untill $2np \geq 4$, then following the same line in \eqref{eq_3_25} can get the desired conclusion.
Then the proof is complete.
\end{proof}
\begin{Rem}
Theorem \ref{theorem3_1} shows that the semi-tamed Milstein approximations \eqref{eq_1_4}
converge with the standard convergence order one to the exact solution of the SDEs \eqref{eq_1_1} in the strong $L^p$ sense.
\end{Rem}

\section{Exponential mean square stability of the semi-tamed Milstein method}

In this section, we study the exponential mean square stability of the semi-tamed Milstein method. The following assumption is needed to get our stability theorem.
\begin{Assu}
	\label{assumption_4_1}
	Assume that $\phi(0) = \varphi(0) = g(0) = 0$, and there exist nonnegative constants $K, \rho, \theta, v, \bar{v}, \alpha > 1$ with $2\rho > \theta^2$ and $2v > \bar{v}$ such that
	\begin{align}
		\label{eq_4_2}\langle x - y, \phi(x) - \phi(y) \rangle &\leq -\rho\|x - y\|^2,\qquad \|\phi(x) - \phi(y)\| \leq K\|x - y\|, \\
		\label{eq_4_3}\langle x, \varphi(x) \rangle &\leq -v\|x \|^{\alpha+1},\qquad \|\varphi(x)\| \leq \bar{v}\|x\|^\alpha,  \\
		\label{eq_4_4}\|g(x) - g(y)\| &\leq \theta\|x - y\|,\qquad \|L^jg_{i}(x) - L^jg_{i}(y)\| \leq \beta\|x - y\|,
	\end{align}
for all $x, y \in \mathbb{R}^d$.
\end{Assu}

\begin{Theo}\cite{Mao1997, zong2014}
	\label{theorem_4_1}
	Let $f(x)$ and $g(x)$ satisfy the local Lipschitz condition. If there exists a positve constant $\gamma$ such that
	\begin{align}
		\label{eq_4_5}
		2\langle x, f(x) \rangle + \|g(x)\|^2 \leq -\gamma\|x\|^2,\qquad \forall x \in \mathbb{R}^d,
	\end{align}
then the solution of Eq. \eqref{eq_1_1} satisfies
	\begin{align}
		\label{eq_4_6}
		\mathbb{E}\|X_t\|^2 \leq \mathbb{E}\|\xi\|^2e^{-\gamma t}.
	\end{align}
\end{Theo}

Obviously, if Assumption \ref{assumption_4_1} holds and $\gamma = 2\rho - \theta^2$, then \eqref{eq_4_6} holds for \eqref{eq_1_1} by Theorem 4.1.

\begin{Theo}
	\label{theorem_4_2}
	Let conditions $\eqref{eq_4_2}-\eqref{eq_4_4}$ be fulfilled. Then there exists a stepsize $h^*$ such that for any $h < h^*$, the numerical approximation \eqref{eq_1_4} has the following property
	\begin{equation}
		\label{eq_4_7}
		\mathbb{E}\|Y^N_n\|^2 \leq \mathbb{E}\|\xi\|^2e^{-nh\gamma_h},
	\end{equation}
where $\gamma_h > 0$ and it satisfies $\lim\limits_{h\rightarrow 0}\gamma_h = 2\rho - \theta^2.$
\end{Theo}

\begin{proof}
Taking square on both sides of
\eqref{eq_1_4} and using the notation $A_n$ given in \eqref{eq_2_12}, we obtain
	\begin{align}
		\label{eq_4_8}
\|Y^N_{n+1}\|^2 &= \|Y^N_n\|^2 + \|\phi(Y^N_n)\|^2h^2 + \frac{\|\varphi(Y^N_n)\|^2h^2}{(1 + \|\varphi(x)\|h)^2} + \|g(Y^N_n)\Delta W_n\|^2 + \|A_n\|^2 + 2h\langle Y^N_n, \phi(Y^N_n)\rangle  \notag\\
		      &\quad + \frac{2\langle Y^N_n + \phi(Y^N_n)h, \varphi(Y^N_n)h\rangle}{1 + \|\varphi(x)\|h} + 2\bigg\langle Y^N_n + \phi(Y^N_n)h + \frac{\varphi(Y^N_n)h}{1 + \|\varphi(x)\|h}, g(Y^N_n)\Delta W_n \bigg\rangle   \notag\\
			  &\quad + 2\bigg\langle Y^N_n + \phi(Y^N_n)h + \frac{\varphi(Y^N_n)h}{1 + \|\varphi(x)\|h}, A_n \bigg\rangle + 2\langle g(Y^N_n)\Delta W_n, A_n \rangle.
	\end{align}
From condition \eqref{eq_4_3}, we have
\begin{equation}
	\label{eq_4_9}
	\frac{2\langle Y^N_n + \phi(Y^N_n)h, \varphi(Y^N_n)h\rangle}{1 + \|\varphi(x)\|h} \leq 2\frac{-v\|Y^N_n\|^{\alpha+1}h}{1 + \|\varphi(x)\|h} + 2\frac{\|\phi(Y^N_n)\|\cdot\|\varphi(Y^N_n)\|h^2}{1 + \|\varphi(x)\|h} \leq 2\frac{-(v - K\bar{v}h)\|Y^N_n\|^{\alpha+1}h}{1 + \|\varphi(x)\|h}.
\end{equation}
Additionally, we define $\bar{\Omega}_{n} = \{ \omega \in \Omega: \|Y^N_n\|\geq 1\}$, then we obtain
\begin{equation}
	\label{eq_4_10}
	\mathds{1}_{\bar{\Omega}_{n}}\frac{\|\varphi(Y^N_n)\|^2h^2}{(1 + \|\varphi(x)\|h)^2} \leq \mathds{1}_{\bar{\Omega}_{n}}\frac{\|\varphi(Y^N_n)\|h}{1 + \|\varphi(x)\|h} \leq \mathds{1}_{\bar{\Omega}_{n}}\frac{\bar{v}\|Y^N_n\|^{\alpha+1}h}{1 + \|\varphi(x)\|h},
\end{equation}
and
\begin{equation}
	\label{eq_4_11}
	\mathds{1}_{\bar{\Omega}_{n}^c}\frac{\|\varphi(Y^N_n)\|^2h^2}{(1 + \|\varphi(x)\|h)^2} \leq \mathds{1}_{\bar{\Omega}_{n}^c}\frac{\bar{v}^2\|Y^N_n\|^{2\alpha}h^2}{1 + \|\varphi(x)\|h} \leq \mathds{1}_{\bar{\Omega}_{n}^c}\frac{\bar{v}^2\|Y^N_n\|^{\alpha+1}h^2}{1 + \|\varphi(x)\|h}.
\end{equation}
Define $h_1 = \frac{2v-\bar{v}}{2Kv} \wedge \frac{2v}{(2K+\bar{v})\bar{v}}$. Then for $h < h_1$, inserting \eqref{eq_4_9}-\eqref{eq_4_11} into \eqref{eq_4_8}, we have
\begin{align}
	\label{eq_4_12}
	\|Y^N_n\|^2 &\leq \|Y^N_n\|^2 + \|A_n\|^2 + \|\phi(Y^N_n)\|^2h^2 + \|g(Y^N_n)\Delta W_n\|^2 + 2h\langle Y^N_n, \phi(Y^N_n)\rangle  \notag\\
	          &\quad + 2\bigg\langle Y^N_n + \phi(Y^N_n)h + \frac{\varphi(Y^N_n)h}{1 + \|\varphi(x)\|h}, g(Y^N_n)\Delta W_n \bigg\rangle  \notag\\
	          &\quad + 2\bigg\langle Y^N_n + \phi(Y^N_n)h + \frac{\varphi(Y^N_n)h}{1 + \|\varphi(x)\|h}, A_n \bigg\rangle  \notag\\
	          &\quad + 2\langle g(Y^N_n)\Delta W_n, A_n \rangle.
\end{align}
Moreover, using mutual independent of $\Delta W_n^{j_1}, \Delta W_n^{j_2}$ and the fact that $Y^N_n$ is independent of $\Delta W_n^{j}$, we have
\begin{align}
	\label{eq_4_13}
	\mathbb{E}&\bigg\langle Y^N_n + \phi(Y^N_n)h + \frac{\varphi(Y^N_n)h}{1 + \|\varphi(x)\|h}, A_n \bigg\rangle \notag\\
	&= \displaystyle\sum_{\substack{j_1,j_2=1\\j_1\neq j_2}}^{m}\mathbb{E}\bigg\langle Y^N_n + \phi(Y^N_n)h + \frac{\varphi(Y^N_n)h}{1 + \|\varphi(x)\|h}, \frac{1}{2}L^{j_1}g_{j_2}(Y^N_n)\Delta W_n^{j_1}\Delta W_n^{j_2} \bigg\rangle  \notag\\
	&+ \displaystyle\sum_{\substack{j=1}}^{m}\mathbb{E}\bigg\langle Y^N_n + \phi(Y^N_n)h + \frac{\varphi(Y^N_n)h}{1 + \|\varphi(x)\|h}, \frac{1}{2}L^{j}g_{j}(Y^N_n)((\Delta W_n^{j})^2 - h) \bigg\rangle  \notag\\
	&= 0.
\end{align}
Similarly, we have
\begin{equation}
	\label{eq_4_14}
	\mathbb{E}\langle g(Y^N_n)\Delta W_n, A_n \rangle = 0,
\end{equation}
and
\begin{align}
	\label{eq_4_15}
	\mathbb{E}\|A_n\|^2 &\leq \frac{m^2}{4}\mathbb{E}\displaystyle\sum_{j_1,j_2=1}^{m}\|L^{j_1}g_{j_2}(Y^N_n)\|^2|\Delta W^{j_1}_n\Delta W^{j_2}_n - \delta_{j_1,j_2}h|^2  \notag\\
			   &\leq \frac{m^2}{4}\displaystyle\sum_{\substack{j_1,j_2=1\\j_1\neq j_2}}^{m}\mathbb{E}\|L^{j_1}g_{j_2}(Y^N_n)\|^2|\Delta W^{j_1}_n\Delta W^{j_2}_n|^2  \notag\\
			   &\quad + \frac{m^2}{4}\displaystyle\sum_{j=1}^{m}\mathbb{E}\|L^{j}g_{j}(Y^N_n)\|^2|\Delta W^{j}_n\Delta W^{j}_n - h|^2  \notag\\
			   &\leq \frac{m^2}{4}(m^2 + m)\beta h^2\mathbb{E}\|Y^N_n\|^2.
\end{align}
Taking expectation on both sides of \eqref{eq_4_12} and using \eqref{eq_4_13}, \eqref{eq_4_14}, \eqref{eq_4_15} as well as Assumption \eqref{assumption_4_1}, we obtain
\begin{align}
	\label{eq_4_16}
	\mathbb{E}\|Y^N_n\|^2 &= \mathbb{E}\|Y^N_n\|^2 - 2\rho h\mathbb{E}\|Y^N_n\|^2 + Kh^2\mathbb{E}\|Y^N_n\|^2 + \theta^2h\mathbb{E}\|Y^N_n\|^2 + \frac{m^2}{4}(m^2 + m)\beta h^2\mathbb{E}\|Y^N_n\|^2  \notag\\
						&= \mathbb{E}\|Y^N_n\|^2\left(1 - 2\rho h + Kh^2 + \theta^2h + \frac{m^2}{4}(m^2 + m)\beta h^2  \right)  \notag\\
						&\leq \mathbb{E}\|Y^N_{n-1}\|^2\left(1 - 2\rho h + Kh^2 + \theta^2h + \frac{m^2}{4}(m^2 + m)\beta h^2  \right)^2  \notag\\
						&\leq \cdots  \notag\\
						&\leq \left(1 - 2\rho h + Kh^2 + \theta^2h + \frac{m^2}{4}(m^2 + m)\beta h^2  \right)^{n+1}\mathbb{E}\|\xi\|^2  \notag\\
						&\leq \mathbb{E}\|\xi\|^2e^{-\gamma_h(n+1)h},
\end{align}
for any $h < h* := h_1 \wedge \frac{2\rho - \theta^2}{\frac{m^2}{4}(m^2 + m)\beta + K},$ where $\gamma_h = 2\rho - \theta^2 + (K + \frac{m^2}{4}(m^2 + m)\beta)h $ and $\gamma_h > 0.$
The proof is complete.
\end{proof}

\begin{Rem}
Theorem \ref{theorem_4_2} shows that the semi-tamed Milstein scheme \eqref{eq_1_4} can reproduce
the exponential mean square stability of the exact solution of the SDEs \eqref{eq_1_1}.
\end{Rem}

\section{Numerical Result}
This section is devoted to presenting numerical results that illustrate the above theoretical analysis. The corresponding Matlab source codes can be found on GitHub \cite{github}.

\begin{Exa}
\begin{equation}
	\label{eq_5_1}
	dX_t = (2X_t - X_{t}^5)dt + X_tdW_t,\qquad X_0 = 1.
\end{equation}
\end{Exa}

\begin{figure}[htb]
	\centering
	\includegraphics[width=0.5\textwidth]{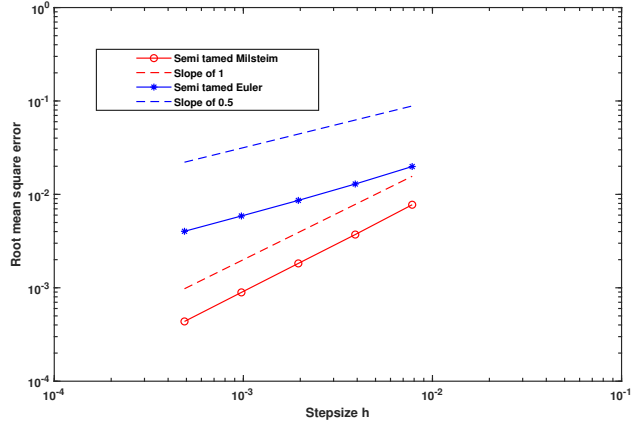}
	\caption{Root mean approximation error versus stepsize $h$ to approxiamtion \eqref{eq_1_4} }
	\label{fig:conv}
\end{figure}
\begin{figure}[htb]
	\centering
	\includegraphics[width=0.6\textwidth]{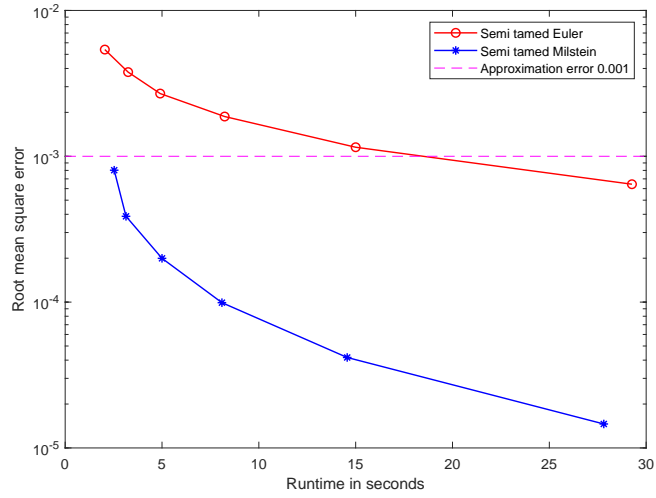}
	\caption{Root mean square approximation error versus runtime for $N\in\{2^{11},\cdots,2^{16}\}$}
	\label{fig_runtime}
\end{figure}
The drift coefficient $f(x) = 2x - x^5$ in \eqref{eq_5_1} does not satisfy the global Lipschitz condition, but it can be divided into a Lipschitz continous term $\phi(x) = 2x$ and a non-Lipschitz continuous term $\varphi(x) = -x^5$. It follows from Theorem \ref{theorem3_1} that the mean-square convergence order of the semi-tamed Milstein method applying to Eq.\eqref{eq_5_1} equals $1$. In Fig.\ref{fig:conv}, the red solid line depicts the root mean-square errors $(\mathbb{E}\left\|X_T - Y^N_N\right\|^2)^{1/2}$ of the semi-tamed Milstein as a function of the stepsize $h$ in log-log plot, while the blue solid line is for the semi-tamed Euler method. The expectation here is approximated by the mean of $5000$ independent realizations. Fig.\ref{fig:conv} shows that the semi-tamed Milstein method gives an error that decreases proportional to $h$ as expected, while the semi-tamed Euler method gives errors that decrease proportional to $h^{1/2}$. We test this by assuming the root mean square approximation errors $e_h$ obey a power-law relation $e_h=C h^r$ for $C>0$ and $r>0$, that is, $\log e_h=\log C+r\log h$. Then we do a least squares power law fit, and obtain $r=1.0288, resid=0.0573$ for the semi-tamed Milstein method, where $resid$ is the least squares residual. The expected convergence rate is confirmed by these results.

We present in Fig.\ref{fig_runtime} the root mean square errors of the semi-tamed Euler method and the semi-tamed Milstein method as a function of the run-time when $N\in\{2^{11},\cdots,2^{16}\}$ to show the efficiency of the semi-tamed Milstein method. Both two methods are calculated on our Quad-Core Intel processor running at 3.5GHz. The expectation here is approximated by the mean of $1000$ independent paths. In Fig.\ref{fig_runtime}, it turns out that $N = 2^{11}$ in the case of the semi-tamed Milstein method and that $N = 2^{16}$ in the case of semi-tamed Euler method achieve the desired precision $\epsilon = 0.001$. Furthermore, the semi-tamed Milstein method requires 2.5313s while the semi-tamed Euler requires 29.2656s to achieve the precision $\epsilon = 0.001$, that is the semi-tamed Milstein method has more advantages in computational efficiency than the semi-tamed Euler method.

Next we investigate the exponential mean square stability of the semi-tamed Milstein method \eqref{eq_1_4}. For simplicity, we consider a one-dimensional SDE.

\begin{Exa}
\begin{equation}
	\label{eq_5_2}
		dX_t = (-2X_t - X_{t}^5)dt + \sqrt{2}X_tdW_t,\qquad X_0 = 1.
\end{equation}
\end{Exa}

From Theorem \ref{theorem_4_1}, we can deduce that the exact solution of \eqref{eq_5_2} is exponentially mean square stable, since
\begin{equation}
	\label{eq_5_3}
	2x^Tf(x) + |g(x)|^2 = -4x^2 - 2x^6 + 2x^2 \leq -2x^2.
\end{equation}

\begin{figure}[htp]
	\centering
	\subfloat[Semi-tamed Euler]{
		\includegraphics[width=0.4\textwidth,height=0.15\textheight]{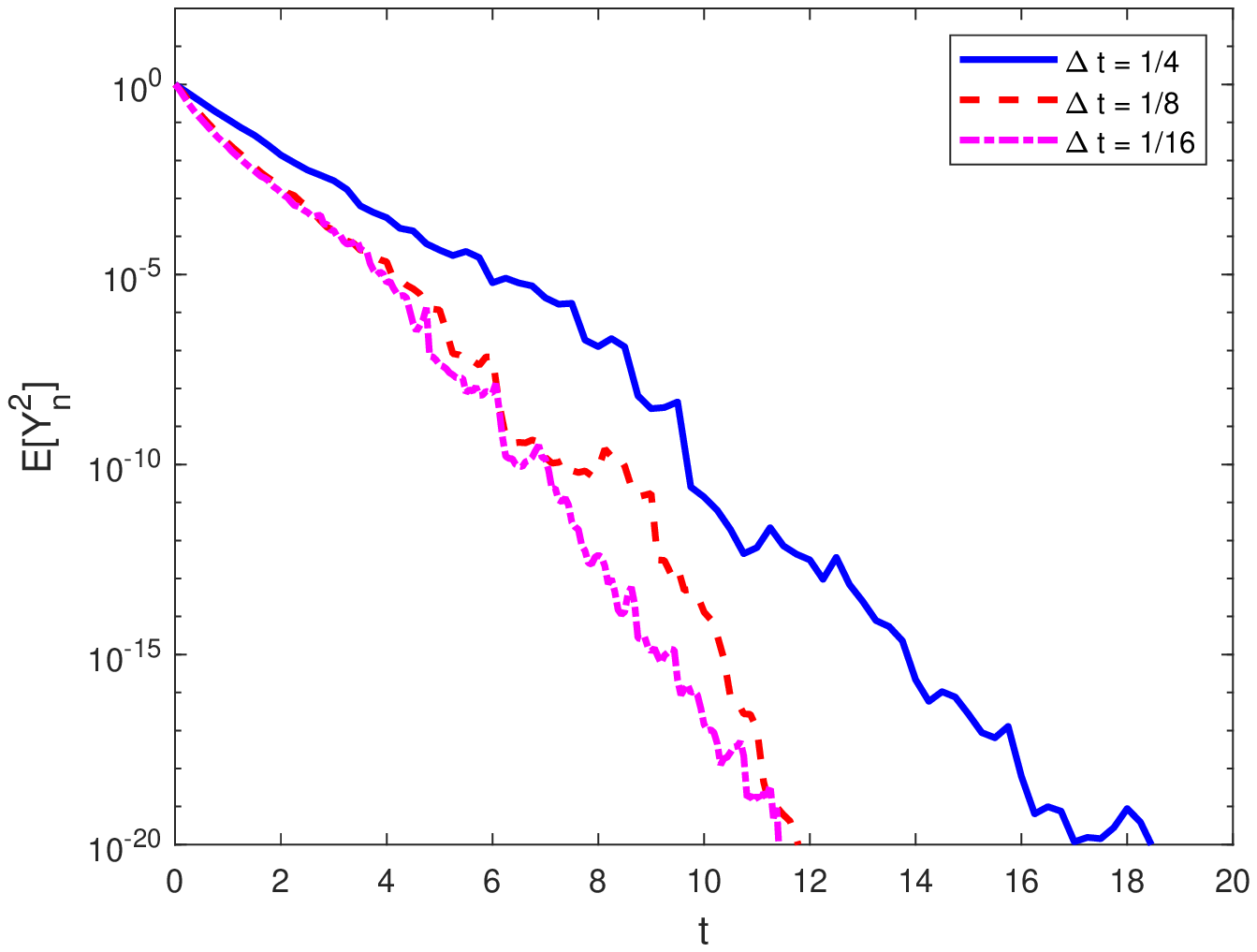}
	}\hspace{.2in}
	\subfloat[Semi-tamed Milstein]{
		\includegraphics[width=0.4\textwidth,height=0.15\textheight]{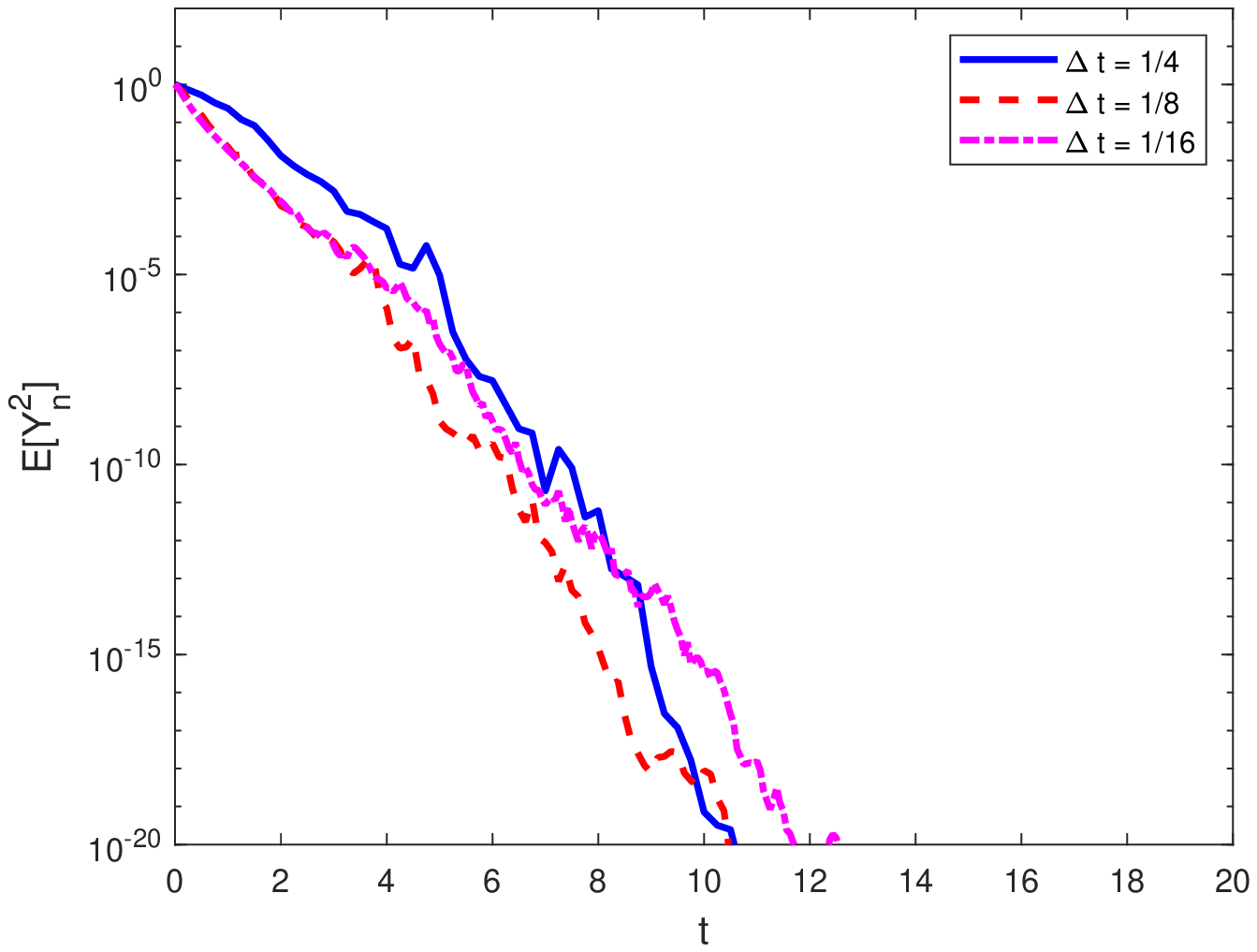}
	}\\
	\subfloat[Tamed Euler]{
		\includegraphics[width=0.4\textwidth,height=0.15\textheight]{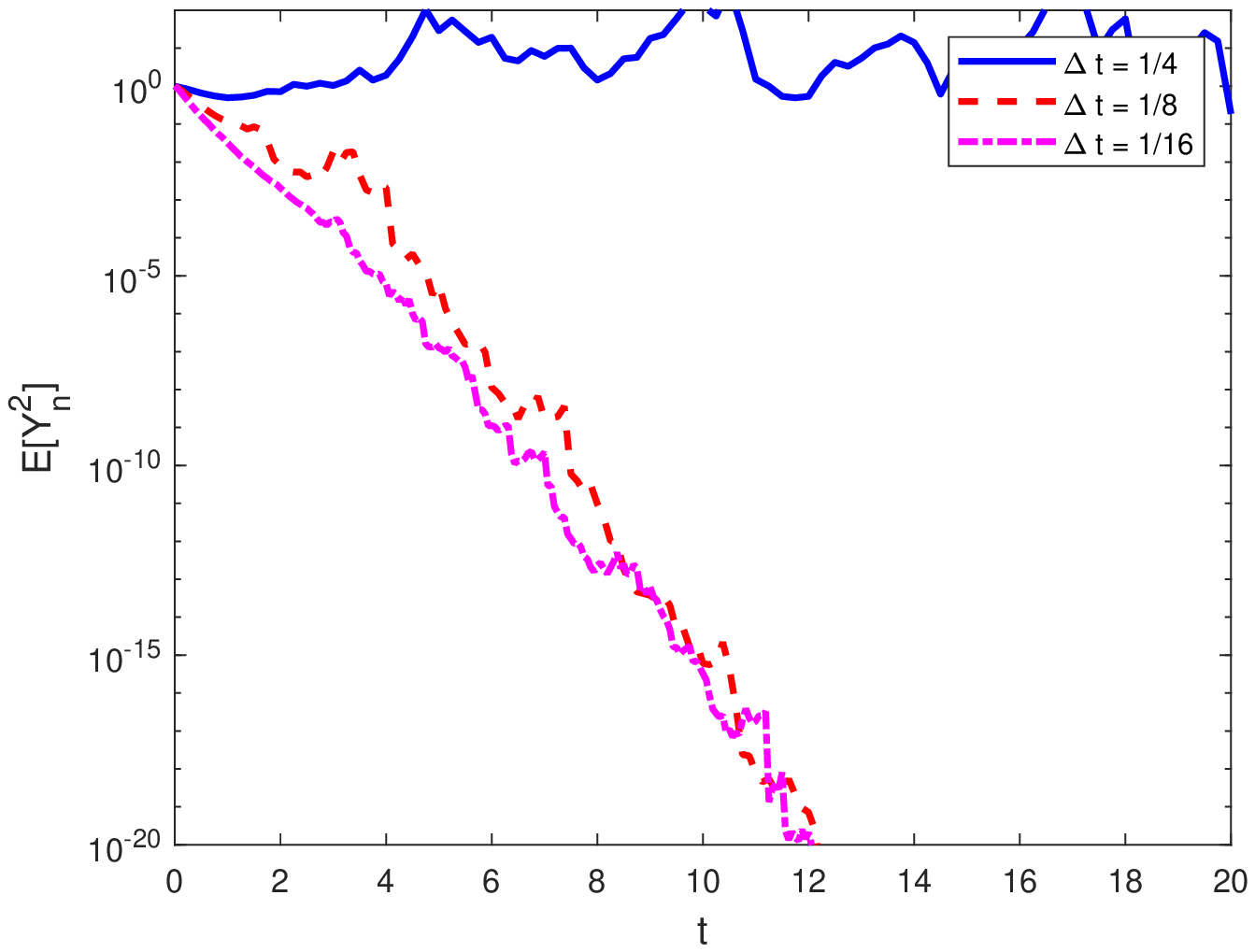}
	}\hspace{.2in}
	\subfloat[Tamed Milstein]{
		\includegraphics[width=0.4\textwidth,height=0.15\textheight]{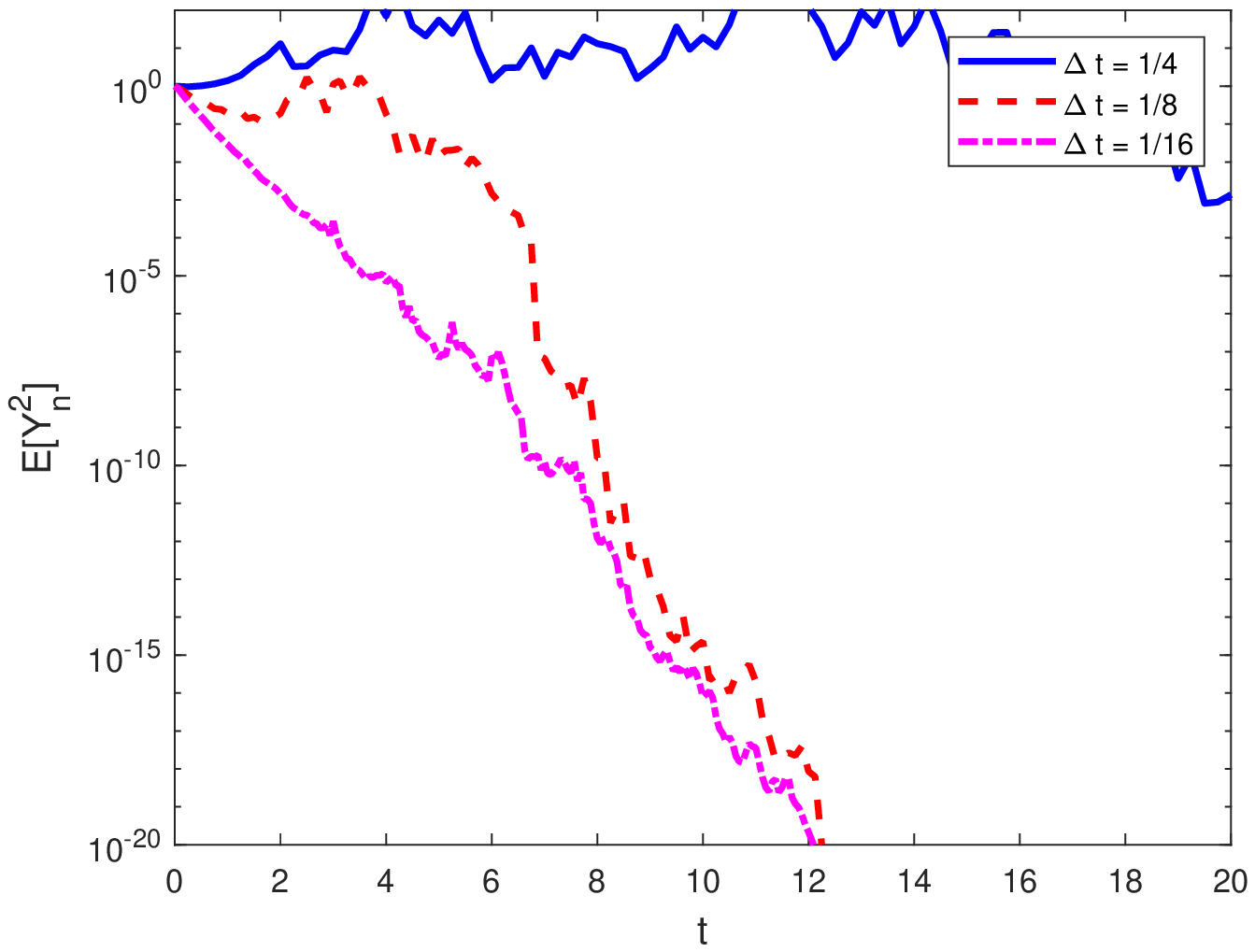}
	}
	\caption{Numerical simulation of $E|X_t^2|$, using different schemes and different stepsizes.}
	\label{fig_stability}
\end{figure}

 In order to test the mean square stability of the tamed Euler, semi-tamed Euler, tamed Milstein and semi-tamed Milstein methods, we take different stepsizes $h = 1/4, 1/8, 1/16$ and generate $5000$ independent sample paths for each numerical scheme. From Fig.\ref{fig_stability}, it's found that the semi-tamed Euler method and semi-tamed Milstein method work better than other two schemes under stepsize $h = 1/4$. But our semi-tamed Mistein method has a higher convergence order than the semi-tamed Euler method as shown in Theorem \ref{theorem3_1} and Fig.\ref{fig:conv}. In addition, the semi-tamed Milstein method is stable with stepsize $h<h^*=0.25$ by Theorem \ref{theorem_4_2}. Fig.\ref{fig_stability} shows that the stability bound obtained in section 4 may be not optimal.

\noindent {\bf Acknowledgments.} The authors would like to thank Prof.Siqing Gan and Prof.Xiaojie Wang for their helpful comments.


\begin{thebibliography}{00}

\bibitem{Mao1997}\label{bibe_16} X. Mao, Stochastic Differential Equations and their Applications, Horwood, Chichester, 1997.

\bibitem{zong2020}\label{bibe_22} X. Zong, T. Li, G. Yin, J. Zhang, Delay tolerance for stable stochastic systems and extensions, IEEE Trans. Automat. Contr. 66(6) (2021) 2604-2619.

\bibitem{maruyama1955}\label{bibe_5} G. Maruyama, Continuous Markov processes and stochastic equations, Rend. Circolo. Math. Palermo. 4(1) (1955) 48-90.

\bibitem{kloeden1992}\label{bibe_4}P. E. Kloeden, E. Platen, Numerical Solution of Stochastic Differential Equations, Springer, Berlin, 1992.

\bibitem{milstein1995}\label{bibe_6} G. N. Milstein, Numerical Integration of Stochastic Differential Equations, Kluwer Academic, Dordrecht, 1995.

\bibitem{Yuan-Mao2008}\label{bibe_21} C. Yuan, X. Mao, A note on the rate of convergence of the Euler-Maruyama method for stochastic differential equations, Stoch. Anal. Appl. 26 (2008) 325-333.

\bibitem{hutzenthaler2009}\label{bibe_7} M. Hutzenthaler, A. Jentzen, Non-globally Lipschitz counter examples for the stochastic Euler scheme, ArXiv preprint arXiv:0905.0273v1, 2009.

\bibitem{hutzenthaler2011}\label{bibe_8} M. Hutzenthaler, A. Jentzen, P.E. Kloeden, Strong and weak divergence in finite time of Euler's method for stochastic differential equations with non-globally Lipschitz continuous coefficients, Pro. R. Soc. Lond. Ser. A Math., Phys. Sci. 467(2130) (2011) 1563-1576.

\bibitem{hutzenthaler2012}\label{bibe_1}M. Hutzenthaler, A. Jentzen, P. E. Kloeden, Strong convergence of an explicit numerical method for SDEs with nonglobally Lipschitz continuous coefficients, Ann. Appl. Probab. 22(4) (2012) 1611-1641.

\bibitem{wang2013}\label{bibe_3} X. Wang, S. Gan, The tamed Milstein method for commutative stochastic differential equations with non-globally Lipschitz continuous coefficients, J. Differ. Equ. Appl. 19(3) (2013) 466-490.

\bibitem{sabanis2016}\label{bibe_14} S. Sabanis, Euler approximations with varying coefficients: the case of superlinearly growing diffusion coefficients, Ann. Appl. Probab. 26(4) (2016) 2083-2105.

\bibitem{mao2015}\label{bibe_11} X. Mao, The truncated Euler-maruyama method for stochastic differential equations, J. Comput. Appl. Math. 290 (2015), 370-384.

\bibitem{guo2018}\label{bibe_12} Q. Guo, W. Liu, X. Mao, R. Yue, The truncated Milstein method for stochastic differential equations with commutative noise, J. Comput. Appl. Math. 338 (2018) 298-310.

\bibitem{zhang2019} W. Zhang, X. Yin, M. Song, M. Liu, Convergence rate of the truncated Milstein method of stochastic differential delay equations, Appl. Math. Comput. 357 (2019) 263-281.



\bibitem{zong2014}\label{bibe_2} X. Zong, F. Wu, C. Huang, Convergence and stability of the semi-tamed Euler scheme for stochastic differential equations with non-Lipschitz continuous coefficients, Appl. Math. Comput. 228 (2014) 240-250.

\bibitem{Gan2014}\label{bibe_20} S. Gan, A. Xiao and D. Wang, Stability of analytical and numerical solutions of nonlinear stochastic delay differential equations, J. Comput. Appl. Math. 268 (2014) 5-22.

\bibitem{Yao2018}\label{bibe_19} J. Yao, S. Gan, Stability of the drift-implicit and double-implicit Milstein schemes for nonlinear SDEs, Appl. Math. Comput. 339 (2018) 294-301.

\bibitem{Burrage1979}\label{bibe_17} K. Burrage, J. C. Butcher, Stability criteria for implicit Runge-Kutta Methods, SIAM J. Numer Anal. 16(1) (1979) 46-57.

\bibitem{gyongy1998}\label{bibe_18} I. Gy{\"o}ngy, A note on Euler's approximations, Potential Anal. 8(3) (1998) 205-216.

\bibitem{hutzenthaler2011b}\label{bibe_9} M. Hutzenthaler, A. Jentzen, Convergence of the stochastic Euler scheme for locally Lipschitz coefficients, Found. Comput. Math. 11(6) (2011) 657-706.

\bibitem{schmallfus1997}\label{bibe_10} B. Schmallfu{\ss}, The random attractor of the stochastic Lorenz system, Z. Angew. Math. Phys. 48(6) (1997) 951-975.

\bibitem{prato1992}\label{bibe_13} G. D. Prato, J. Zabczyk, Stochastic Equations in Infinite Dimensions, Cambridge university press, Cambridge, 1992.

\bibitem{github}\label{bibe_23} $https://github.com/xiaohanWan/issue\_Matlab\_Program$.

\end{thebibliography}
\end{document}